\documentclass[reqno]{amsart}
\usepackage{amsmath,amsthm,amssymb,amsfonts,mathrsfs,amscd}
\usepackage{multirow}
\usepackage{graphicx}
\usepackage{bm}
\usepackage{mathrsfs}
\usepackage{dsfont}
\usepackage{tikz}
\usepackage{tikz-3dplot}
\usepackage{subfig}
\usepackage{stmaryrd}

\usepackage{float}
\usepackage[textfont=it]{caption}%[justification=raggedright]

\newtheorem{remark}{Remark}[section]

\newtheorem{lemma}{Lemma}[section]
\newtheorem{theorem}{Theorem}[section]

\makeatletter
  \@addtoreset{equation}{section}
  \newcommand\figcaption{\def\@captype{figure}\caption}
  \newcommand\tabcaption{\def\@captype{table}\caption}
\makeatletter

\begin{document}

\title{A novel least squares method for Helmholtz equations with large wave numbers}
%in inhomogeneous media

\author{Qiya Hu}
\author{Rongrong Song}

\thanks{1. LSEC, ICMSEC, Academy of Mathematics and Systems Science, Chinese Academy of Sciences, Beijing
100190, China; 2. School of Mathematical Sciences, University of Chinese Academy of Sciences, Beijing 100049,
China (hqy@lsec.cc.ac.cn, songrongrong@lsec.cc.ac.cn). This work was funded by Natural Science Foundation of China G11571352.}

\maketitle

{\bf Abstract.}
In this paper we are concerned with numerical methods for Helmholtz equations with large wave numbers. We design
a least squares method for discretization of the considered Helmholtz equations. In this method, an auxiliary unknown is introduced on the common
interface of any two neighboring elements and a quadratic objective functional is defined by the jumps
of the traces of the solutions of local Helmholtz equations across all the common interfaces, where the local Helmholtz equations are defined on elements
and are imposed Robin-type boundary conditions given by the auxiliary unknowns. A minimization problem with the objective functional is proposed
to determine the auxiliary unknowns. The resulting discrete system of the auxiliary unknowns is Hermitian positive definite and so it can be solved by the preconditioned conjugate
gradient (PCG) method. Under some assumptions we show that the generated approximate solutions possess almost the same $L^2$ convergence order as the plane wave methods
(for the case of constant wave number).
Moreover, we construct a substructuring preconditioner for the discrete system of the auxiliary unknowns. Numerical experiments show that the proposed methods
are very effective and have little ``wave number pollution" for the tested Helmholtz equations with large wave numbers.



\pagestyle{myheadings}
\thispagestyle{plain}
\markboth{}{QIYA HU AND RONGRONG SONG}

\section{Introduction}
%For simplicity of exposition, we only consider two dimensional problem in this paper.
Let $\Omega$ be a bounded, connected and Lipschitz domain in $\mathbb{R}^2$. Consider the Helmholtz equations
\begin{equation} \label{eq}
\left\{
\begin{aligned}
& -\Delta u - {\mathbf \kappa}^2 u = f
& \text{in}\ \Omega,\\
& \frac{\partial u}{\partial {\bf n}} + i{\mathbf \kappa} u = g
& \text{on}\ \partial\Omega,
\end{aligned}
\right.
\end{equation}
where ${\bf n}$ denotes the unit outward normal on the boundary $\partial\Omega$ and ${\mathbf \kappa}$ is the wave number defined by ${\bf \kappa}({\bf x})=\frac{\omega}{c({\bf x})}>0$, with $\omega>0$ being a constant and $c({\bf x})$ being a bounded and positive function defined on $\Omega$. In applications, $\omega$ denotes the angular frequency, which may be very large, and $c({\bf x})$ denotes the wave speed (the acoustic velocity),
which may not be a constant function on $\Omega$, i.e., the involved media is inhomogeneous.
%The above Robin boundary conditions can be replaced by other more complicated boundary conditions, with some obvious modifications in the methods proposed in this paper.

Helmholtz equation is the basic model in sound propagation. It is a very important
topic to design a high accuracy method for Helmholtz equations with large wave numbers,
such that the so called ``wave number pollution" can be reduced. The ``wave number pollution" says that, for a finite element method for the discretization of (\ref{eq}), the mesh size $h$ must satisfy $h\omega^{1+\delta}=O(1)$ for some positive number $\delta$ to achieve a given accuracy of the approximate solutions when the wave number $\omega$ increases, which means that the
accuracies of the approximate solutions are obviously destroyed if fixing the value of $h\omega$ but increasing the wave number $\omega$. For convenience, we call the parameter $\delta$ as the ``pollution index", which describes
the degree of wave number pollution. For the standard linear finite element method, the ``pollution index" $\delta=1$ (see \cite{fem}). Of curse, we hope to design a ``good" finite element method for (\ref{eq}) such that the pollution index $\delta$ is sufficiently small.

In recent years, many interesting methods for the discretization of Helmholtz equations with large wave numbers have been proposed, for example (but not all), the higher order finite element methods (hp-FEM) \cite{DuWu2015,fem}, the ultra weak variational formulation (UWVF) \cite{ref11}, the plane wave least squares (PWLS) methods \cite{HuY2014,Hu2014, Monk1999}, the plane wave discontinuous Galerkin (PWDG) methods \cite{Gittelson2009, Hiptmair2011},
the method of fundamental solutions \cite{Alves2005-2, Deckers2014}, the plane wave method with Lagrange multipliers (PWLM) \cite{Farhat2003}, the variational theory of complex rays \cite{Riou2008},
the high order element discontinuous Galerkin method (HODG) \cite{DuWu2015,Wu2011}, local discontinuous Galerkin method (LDG) \cite{Feng2013}, hybridizable discontinuous Galerkin method (HDG) \cite{Chen2013, Cui2013, Gopalakrishnan2015, Gries2011, Monk2010, Nguyena2015, Stang2016}
and the discontinuous Petrov-Galerkin (DPG) method \cite{Demk2012, Gopalakrishnan2014, Yuan2015}, the ray-based finite element method \cite{Fang2016} and the generalized plane wave method \cite{IMbMon2015}.
All these methods are superior to the standard linear finite element method in the sense that the pollution index $\delta<1$.

It is known that the plane wave finite element methods have little ``wave number pollution" (i.e., the pollution index $\delta$ is very small) and can generate higher accuracy approximations than the polynomial basis finite element methods for solving the Helmholtz equations with large (piecewise constant) wave numbers when finite element spaces have the same degrees of freedom. A comparison of finite element methods based on high-order polynomial basis
functions and plane wave basis functions was given in \cite{Lieus2016}. The numerical results reported in \cite{Lieus2016} indicate that, if only the degrees of freedom on element boundaries for high-order polynomial method are calculated (the degrees of freedom in the interior of elements are eliminated), the high-order polynomial method can deliver comparable to the PWDG method.
Unfortunately, the plane wave methods cannot be directly applied to the discretization of nonhomogeneous Helmholtz equations in inhomogeneous media.
A plane wave method combined with local spectral element for nonhomogeneous Helmholtz equations in homogeneous media was proposed in \cite{HuY2017} (see also \cite{Hu2014}).
A generalized plane wave method for homogeneous Helmholtz equations in inhomogeneous media was introduced in \cite{IMbMon2015}.

The HDG-type methods (and the DPG method) have been studied in many works (see the references listed above). We would like to simply recall the ideas of the HDG methods. Let $\Omega$ be decomposed into a union of elements $\{\Omega_k\}$, and let $\gamma$ denote the element interface, which is a union of all the common edges of two neighboring elements. For the HDG-type methods, the equation (\ref{eq}) is first transformed into a first-order system of the original unknown $u$ and an auxiliary unknown $\Phi=(i\omega)^{-1}\nabla u$, then the restrictions of the unknowns $u$ and $\Phi$ on the elements $\{\Omega_k\}$ are eliminated by solving all the local first-order systems to obtain an interface equation of the trace $u|_{\gamma}$ (and the trace $(\nabla u\cdot {\bf n})|_{\gamma}$ in \cite{Monk2010}) in some manner.
%For the HDG method, the interface unknown is just the trace $u|_{\gamma}$ and the interface equation is still indefinite.
For the DPG method, there are two interface unknowns that are defined by the traces $u|_{\gamma}$ and $(\nabla u\cdot {\bf n})|_{\gamma}$ and the interface equation becomes Hermitian positive definite by introducing
nonstandard test space that is the image of the trial space under a suitable mapping. In both the HDG-type methods and the DPG method, the unknown needed to be globally solved was defined on the interface $\gamma$,
so these methods have less cost of calculation than the standard $hp$ finite element method proposed in \cite{fem}. The HDG-type methods and the DPG method have their respective merits: the HDG-type methods are easier to implement than the DPG method since the HDG-type methods use the standard polynomial basis functions; the interface equation needed to be solved globally is Hermitian positive definite for the DPG method, but it is still indefinite as the original equation (\ref{eq}) for the HDG-type methods.

In the present paper, we design a novel discretization method for Helmholtz equations with large wave numbers such that the method can absorb the merits of the HDG-type methods and the DPG method.
%has little ``wave number pollution" and possesses some other nice features.
The basic ideas of the new method can be roughly described as follows. We introduce an auxiliary unknown $\lambda_h$ that is an edge-wise $q$ order polynomial on $\gamma$,
and compute $p$ order ($p\geq q+2$) polynomial solutions $\{u_{h,k}\}$ of the discrete variational problems of all local Helmholtz equations, where
each local Helmholtz equation is the restriction of (\ref{eq}) on some element $\Omega_k$ and is imposed a Robin-type boundary condition given by the auxiliary unknown $\lambda_h$.
We define a minimization problem with a quadratic objective functional defined by the jumps of the traces of the solutions $\{u_{h,k}\}$ across the interface $\gamma$.
This minimization problem results in a Hermitian positive definite algebraic system of the auxiliary unknown $\lambda_h$. After solving the algebraic system,
we can easily obtain an approximate solution of the original Helmholtz equation by solving small local problems on the elements in parallel manner.
This method has some similarity with the HDG method but it has essential differences from the HDG method:
(a) each element subproblem is just the local variational problem of the original Helmholtz equation (\ref{eq}), so only one internal unknown $u_{h,k}$ needs to be computed for an element $\Omega_k$;
%, however two unknown need to be computed on each element for the HDG method (and the DPG method);
(b) the interface unknown $\lambda_h$, which may be discontinuous on the interface $\gamma$, is defined independently on every edge of elements; (c) the interface unknown $\lambda_h$
is determined by a minimization problem, so the interface equation is Hermitian positive definite.
%however the interface unknown in the HDG method is defined by the standard continuity condition on element interface.

The new method possesses the following merits: (i) the proposed method is practical to general nonhomogeneous Helmholtz equations in inhomogeneous media (comparing
the plane wave methods); (ii) the algebraic system of $\lambda_h$ is Hermitian positive definite (comparing the HDG-type methods, the PWDG method and the PWLM method),
so it can be solved by the PCG method, which has stable convergence and less cost of calculation, and the construction of preconditioner for this system has more choices (for example,
the well-known BDDC method can be considered); (iii) it is cheap to implement since only one unknown $u_{h,k}$ is introduced in an element $\Omega_k$ and only one unknown $\lambda_h|_{\gamma_{lj}}$
is involved on each local interface $\gamma_{lj}$
%, which means that the unknown $\lambda_h$ needed to be globally solved has less degrees of freedom
(comparing the HDG-type methods and the DPG method); (iv) the method is easy to implement since the subproblem for computing $u_{h,k}$ on an element $\Omega_k$ is directly defined by the original second order Helmholtz equation and the basis
functions on every element $\Omega_k$ and every element edge are standard polynomials (comparing the DPG method).

Since the resulting approximate solution $(u_h,\lambda_h)$ do not satisfy a mixed variational problem (as in the Lagrange multiplier methods) or a hybridizable variational problem
(as in the HDG methods), well-posedness and convergence of the proposed method cannot be proved by the techniques developed in existing works.

By developing some new techniques, we show that the proposed discretization method is well-posed and the resulting approximate solution possesses almost the same $L^2$ error estimate as the
plane wave methods under suitable assumptions, which indicate that the proposed method
has little ``wave number pollution''. In addition, we construct a domain decomposition preconditioner for the algebraic system of $\lambda_h$. The BDDC method is a popular substructuring
domain decomposition method, which was first proposed in \cite{bddc} and then was extended to various models by many researchers. The key idea of the BDDC method is to compute basis functions
of the coarse space by solving local minimization problems. This method has some advantages over the traditional substructuring methods, but
the minimization problems for computing coarse basis functions can be defined only for symmetric and positive definite systems. Thanks to the Hermitian positive definiteness of the algebraic system of $\lambda_h$, we can construct a substructuring preconditioner for the system by the BDDC method. However, we find that the coarse space defined by the BDDC method is unsatisfactory for the current situation. Because of this, we construct a variant of the BDDC preconditioner for the algebraic system of $\lambda_h$ by changing the definition of coarse space. Numerical results indicate that the proposed discretization method and preconditioner are very efficient for the tested Helmholtz equations with large wave numbers.

The paper is organized as follows:
In Section 2, we describe the proposed least squares variational formulation for Helmholtz equations.
In Section 3, we construct a substructuring preconditioner for the discrete system.
The main results about error estimates are presented in Section 4.
In Section 5, we give proofs of the main results in details.
Finally, we report some numerical results to confirm the effectiveness of the new method in Section 6.

\section{A least squares variational formulation}
%We suppose there are $\kappa_0>0$ and $\kappa_1>0$, such that $\kappa_0 \leq \kappa \leq \kappa_1$.

\subsection{Notations}
As usual we partition $\Omega$ into elements in the sense that
$$ \overline{\Omega}=\bigcup_{k=1}^N\overline{\Omega}_k,\quad \Omega_k\bigcap\Omega_j=\emptyset, \quad\text{for } k\not=j. $$
Here each $\Omega_k$ may be curve polyhedron. We use $h_k$ to denote the diameter of $\Omega_k$ and set $h=\max\{h_k\}$.
Let ${\mathcal T}_h$ denote the partition comprised of elements $\{\Omega_k\}_{k=1}^N$. As usual we assume that the partition ${\mathcal T}_h$
is quasi-uniform and regular.

Let $\gamma_{kj}$ denote the common edge of two neighboring elements $\Omega_k$ and $\Omega_j$, and set $\gamma_k=\partial\Omega_k\cap \partial\Omega$
%   = \partial\Omega_k\bigcap\partial\Omega_j \quad\text{for } k\not=j; \quad  \gamma_k=\partial\Omega_k \bigcap \partial\Omega $$
when the intersection is an edge of the element $\Omega_k$. For convenience, define $\gamma = \cup_{k\not=j}\gamma_{kj}$.

Let $q\geq1$ be an integer and choose $p\geq q+2$. Throughout this paper we use the following notations:

% $\bullet\ $ $K$ is a rectangular-shaped element of ${\mathcal T}_h$ and $\partial K$ is its boundary.

% $\bullet\ $ $ V(\Omega) = \prod\limits_{k=1}^N H^1(\Omega_k). $
% = \{ v\in L^2(\Omega):\forall \,K\in{\mathcal T}_h, v_k\in H^1(K) \}

$\bullet\ $ $ V_h^p(\Omega_k) = \{ v\in H^1(\Omega_k): v \text{ is a polynomial whose order does not exceed } p \}. $
% \text{ with respect to each variable}

$\bullet\ $ $ V_h^p({\mathcal T}_h) = \prod_{k=1}^N V_h^p(\Omega_k). $

$\bullet\ $ $ V_h^p(\partial\Omega_k) = \{ v|_{\partial\Omega_k}: v \in V_h^p(\Omega_k) \}. $

$\bullet\ $ $ W(\gamma) = \prod_{k\not=j}^N H^{-\frac{1}{2}}(\gamma_{kj}). $

$\bullet\ $ $ W_h^q(\gamma_{kj}) = \{ \mu\in H^1(\gamma_{kj}):  \mu \text{ is a polynomial whose order does not exceed } q \}. $

$\bullet\ $ $ W_h^q(\gamma) = \prod_{k\not=j} W_h^q(\gamma_{kj}). $

$\bullet\ $ $ W_h^q(\partial\Omega_k \backslash\partial\Omega) = \{ \mu|_{\partial\Omega_k \backslash\partial\Omega}: \mu \in W_h^q(\gamma) \}. $

%$\bullet\ $ $Q_h^p$ denotes the $L^2$ projection operator: $ L^2(\partial\Omega_k) \longmapsto V_h^p(\partial\Omega_k). $

%$\bullet\ $ $Q_h^q$ denotes the $L^2$ projection operator: $ L^2(\partial\Omega_k \backslash\partial\Omega) \longmapsto W_h^q(\partial\Omega_k \backslash\partial\Omega). $

%$\bullet\ $ For a domain $D\subset \Omega$, let $||\cdot||_{1,\kappa,D}$ be defined by
%$$ ||u||_{1,\kappa,D} = \Big( ||\nabla u||_{0,D}^2 + \kappa^2||u||_{0,D}^2 \Big)^\frac{1}{2}. $$

$\bullet\ $ The jump of $v$ across $\gamma_{kj}$: $[v] = v_k-v_j$, where $v$ \text{ is a piecewise smooth function on } ${\mathcal T}_h$ and $v_k=v|_{\Omega_k}$.

$\bullet\ $ $(u, v)_{\Omega_k}=\int_{\Omega_k} u\cdot v ~dx$, \quad $\langle u, v \rangle_{\partial\Omega_k}=\int_{\partial\Omega_k} u\cdot v ~ds.$

%$\bullet\ $ $||u||_{0,\Omega}=\big( \sum\limits_{k=1}^N ||u||^2_{0,\Omega_k} \big)^\frac{1}{2}=\big( \sum\limits_{k=1}^N\int_{\Omega_k} |u|^2 \, dx\big)^\frac{1}{2}$.

%$\bullet\ $ $|u|_{r+1,\Omega}=\big( \sum\limits_{k=1}^N | u|^2_{r+1,\Omega_k}\big)^\frac{1}{2}.$

\subsection{A continuous variational formulation} Let $u\in H^1(\Omega)$.
For each element $\Omega_k$, set $u|_{\Omega_k} = u_k$. For $k>j$, define $\lambda\in W(\gamma)$ as
$$\lambda|_{\gamma_{kj}}
= ( \frac{\partial u_k}{\partial {\bf n}_k}+i\rho u_k )|_{\gamma_{kj}}
= ( -\frac{\partial u_j}{\partial {\bf n}_j}+i\rho u_j )|_{\gamma_{kj}}, $$
where $\rho>0$, ${\bf n}_k$ and ${\bf n}_j$ separately denote the unit outward normal on $\partial\Omega_k$ and $\partial\Omega_j$.
% and we suppose $\rho\ll 1$.
It is clear that the solution $u$ of (\ref{eq}) satisfies the local Helmholtz equation on each element $\Omega_k$ ($k=1,\cdots,N$)
\begin{equation} \label{localeq}
\left\{
\begin{aligned}
& -\Delta u_k - \kappa^2 u_k = f
& \text{in}\ \Omega_k,\\
& \frac{\partial u_k}{\partial {\bf n}_k} \pm i\rho u_k =\pm \lambda
& \text{on}\ \partial\Omega_k\backslash\partial\Omega,\\
& \frac{\partial u_k}{\partial {\bf n}_k} + i\kappa u_k = g
& \text{on}\ \partial\Omega_k\cap\partial\Omega.
\end{aligned}
\right.
\end{equation}
Here the sign $``\pm"$ means that two inverse signs are used on the two side of each local interface $\gamma_{kj}=\partial\Omega_k\cap\partial\Omega_j$:
it takes ``+" on $\gamma_{kj}\subset\partial\Omega_k$, and it takes ``-" on $\gamma_{kj}\subset\partial\Omega_j$.

For each element $\Omega_k$, define the local sesquilinear form
\begin{equation*}
\begin{split}
a^{(k)}(v,w)
& =(\nabla v,\nabla \overline{w})_{\Omega_k}-(\kappa^2 v, \overline{w})_{\Omega_k} \pm i\rho\langle v, \overline{w}\rangle_{\partial\Omega_k \backslash\partial\Omega}\\
& +i\langle \kappa v, \overline{w}\rangle_{\partial\Omega_k \cap\partial\Omega}, \quad v,w\in H^1(\Omega_k)
\end{split}
\end{equation*}
and the local functional
$$ L^{(k)}(v)=(f, \overline{v})_{\Omega_k}+\langle g, \overline{v}\rangle_{\partial\Omega_k \cap\partial\Omega},\quad v\in H^1(\Omega_k).$$
It is easy to see that the variational formulation of (\ref{localeq}) is: to find $u_k(\lambda)\in H^1(\Omega_k)$ such that
\begin{equation}  \label{viar1}
a^{(k)}(u_k(\lambda),v)
=L^{(k)}(v)+\langle\pm \lambda, \overline{v}\rangle_{\partial\Omega_k \backslash\partial\Omega}, \quad\forall\,v\in H^1(\Omega_k).
\end{equation}

We define the quadratic functional
\begin{equation}  \label{minfunc}
J(\mu) = \sum_{\gamma_{kj}} \int_{\gamma_{kj}} |u_k(\mu)-u_j(\mu)|^2 ds, \quad \mu\in W(\gamma)
\end{equation}
and consider the following minimization problem:
find $\lambda \in W(\gamma)$ such that
\begin{equation}
J(\lambda)=\min\limits_{\mu\in W(\gamma)} J(\mu). \label{min}
\end{equation}
It is clear that $u$ is the solution of (\ref{eq}) if and only if $J(\lambda)=0$, which means that $\lambda$ is the solution of the minimization problem (\ref{min}).

In order to give a variational problem of (\ref{min}), we write the solution of (\ref{localeq}) as $u_k(\lambda)= u_k^{(1)}(\lambda)+u_k^{(2)}$, which respectively satisfy
$$ a^{(k)}(u^{(1)}_k(\lambda),v)=\pm\langle \lambda, \overline{v}\rangle_{\partial\Omega_k \backslash\partial\Omega}, \quad\forall\,v\in H^1(\Omega_k) $$
and
$$ a^{(k)}(u^{(2)}_k,v)=L^{(k)}(v), \quad\forall\,v\in H^1(\Omega_k). $$
Then $J(\mu)$ can be written as
$$ J(\mu)=\sum_{\gamma_{kj}} \int_{\gamma_{kj}} |\big(u^{(1)}_k(\mu)-u^{(1)}_j(\mu)\big)+(u^{(2)}_k-u^{(2)}_j)|^2 ds. $$

Define the sesquilinear form
\begin{equation*}
s(\lambda,\mu)=\sum_{\gamma_{kj}} \int_{\gamma_{kj}} (u_k^{(1)}(\lambda)-u_j^{(1)}(\lambda)) \cdot \overline{(u^{(1)}_k(\mu)-u^{(1)}_j(\mu))} ds,\quad \lambda,\mu\in W(\gamma)
\end{equation*}
and the functional
\begin{equation*}
l(\mu)=-\sum_{\gamma_{kj}} \int_{\gamma_{kj}} (u_k^{(2)}-u_j^{(2)}) \cdot \overline{(u^{(1)}_k(\mu)-u^{(1)}_j(\mu))} ds,\quad \mu\in W(\gamma).
\end{equation*}
Therefore the variational problem of the minimization problem (\ref{min}) can be expressed as follows: find $\lambda\in W(\gamma)$ such that
\begin{equation}
s(\lambda,\mu) = l(\mu), \quad \forall ~\mu\in W(\gamma).\label{contin}
\end{equation}

\subsection{The discrete variational formulation} \label{uh}
Let $\lambda_h\in W_h^q(\gamma)$. For each element $\Omega_k$, define $u_{h,k}(\lambda_h)\in V_h^p(\Omega_k)$ by
\begin{equation}  \label{viar2}
a^{(k)}(u_{h,k}(\lambda_h),v_{h})
=L^{(k)}(v_{h})+\langle \pm\lambda_h, \overline{v}_{h}\rangle_{\partial\Omega_k \backslash\partial\Omega}, \quad\forall\,v_{h}\in V_h^p(\Omega_k).
\end{equation}
It is easy to see that the above problem is uniquely solvable.

As in the continuous situation, we decompose $u_{h,k}$ into $u_{h,k} = u_{h,k}^{(1)}(\lambda_h)+u_{h,k}^{(2)}$, which are respectively defined by
\begin{equation*}
a^{(k)}(u_{h,k}^{(1)}(\lambda_h),v_{h})
=\pm\langle \lambda_h, \overline{v}_{h}\rangle_{\partial\Omega_k \backslash\partial\Omega}, \quad\forall\,v_{h}\in V_h^p(\Omega_k)
\end{equation*}
and
\begin{equation*}
a^{(k)}(u_{h,k}^{(2)},v_{h})
=L^{(k)}(v_{h}), \quad\forall\,v_{h}\in V_h^p(\Omega_k).
\end{equation*}
From the computational point of view, the function $u^{(2)}_{h,k}$ can be preliminarily calculated, which will be appeared in the right side of the discrete system,
but the function $u^{(1)}_{h,k}$ cannot be calculated until the function $\lambda_h$ is obtained.

Define the discrete sesquilinear form
\begin{equation*}
s_h(\lambda_h,\mu_h)\!=\!\sum_{\gamma_{kj}}\! \int_{\gamma_{kj}} \! (u_{h,k}^{(1)}(\lambda_h)\!-\!u_{h,j}^{(1)}(\lambda_h)) \!\cdot\! \overline{(u_{h,k}^{(1)}(\mu_h)\!-\!u_{h,j}^{(1)}(\mu_h))} ds, ~~\lambda_h,\mu_h\in W_h^q(\gamma)
\end{equation*}
and the functional
\begin{equation*}
l_h(\mu_h)=-\sum_{\gamma_{kj}} \int_{\gamma_{kj}} (u_{h,k}^{(2)}-u_{h,j}^{(2)}) \cdot \overline{(u_{h,k}^{(1)}(\mu_h)-u_{h,j}^{(1)}(\mu_h))} ds,\quad \mu_h\in W_h^q(\gamma).
\end{equation*}

Therefore the discrete variational problem of (\ref{contin}) can be written as follows: find $\lambda_h\in W_h^q(\gamma)$ such that
\begin{equation}  \label{discre}
s_h(\lambda_h,\mu_h) = l_h(\mu_h), \quad \forall ~\mu_h\in W_h^q(\gamma).
\end{equation}

After $\lambda_h$ is solved from (\ref{discre}), we can easily compute $u_{h,k}$ in parallel by (\ref{viar2}) for every $\Omega_k$. Define
$u_h\in V_h^p(\mathcal{T}_h)$ by $u_h|_{\Omega_k}=u_{h,k}(\lambda_h)$ ($k=1,\cdots,N$). Then $u_h$ should be an approximate solution
of $u$. We would like to emphasize the discrete system (\ref{discre}) has relatively less degrees of freedom, so it is cheaper to be solved.

Let ${\mathcal S}$ be the stiffness matrix associated with the sesquilinear form $s_h(\cdot,\cdot)$, and let $b$ denote the vector associated with $l_h(\cdot)$.
Then the discretization problem (\ref{discre}) leads to the algebraic system
\begin{equation} \label{axb}
{\mathcal S}X = b.
\end{equation}
From the definition of the sesquilinear form $s_h(\cdot,\cdot)$, we know that the matrix ${\mathcal S}$ is Hermitian positive definite, so the system (\ref{axb}) can be solved
by the preconditioned CG method with a positive definite preconditioner. The construction of an efficient preconditioner for ${\mathcal S}$ is an important task (see the next section).
\begin{remark} Since each local finite element space $V_h^p(\Omega_k)$ consists of the standard polynomials, instead of solutions of homogeneous Helmholtz equation in the plane wave methods,
from the viewpoint of algorithm the proposed method is practical to general nonhomogeneous Helmholtz equations in inhomogeneous media.
\end{remark}

\begin{remark} As in the traditional Lagrange multiplier method, we can derive another discrete system of $\lambda_h$ by
the constraints (for all element interfaces $\gamma_{kj}$)
$$ \langle u_{h,k}-u_{h,j},\mu\rangle_{\gamma_{kj}}=0,\quad\forall\mu\in W_h^q(\gamma). $$
However, the coefficient matrix of the resulting system is still indefinite as (\ref{eq}) (comparing the system (\ref{axb})), which makes the solution of the system to be more difficult.
\end{remark}
%In the next section we construct an efficient preconditioner ${\mathcal K}$ for the matrix ${\mathcal S}$.

\section{A domain decomposition preconditioner} \label{pre}

In this section, we are devoted to the construction of a preconditioner ${\mathcal K}$ for ${\mathcal S}$. Thanks to the Hermitian positive definiteness of the matrix ${\mathcal S}$,
we can construct a (Hermitian positive definite) substructuring preconditioner absorbing some ideas in
%by the method proposed in \cite{HuZ2016} and
the BDDC method first introduced in \cite{bddc} (see Section 1 for simple descriptions of the BDDC method). As we will see,
%for the current situation,
the preconditioner designed in this section has essential differences from the one defined in the standard BDDC method.

For convenience, we will define the preconditioner in operator form. To this end, let $S: W_h^q(\gamma)\rightarrow
W_h^q(\gamma)$ denote the discrete operator corresponding to the stiffness matrix ${\mathcal S}$, i.e.,
$$ \langle S\lambda_h,\mu_h\rangle=s_h(\lambda_h,\mu_h),\quad\forall\lambda_h,\mu_h\in W_h^q(\gamma).$$

\begin{figure}[!ht]
\centering
\includegraphics[width=5.3cm,height=5cm]{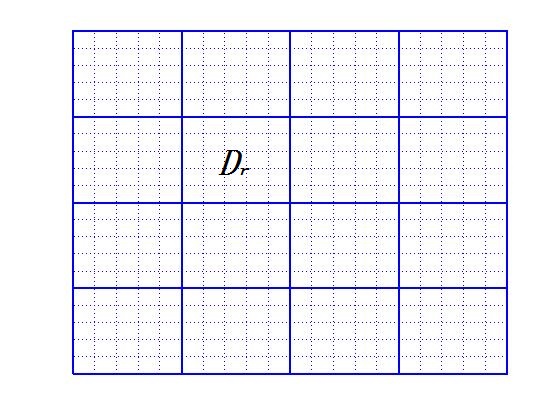}\quad\quad\includegraphics[width=5.3cm,height=5cm]{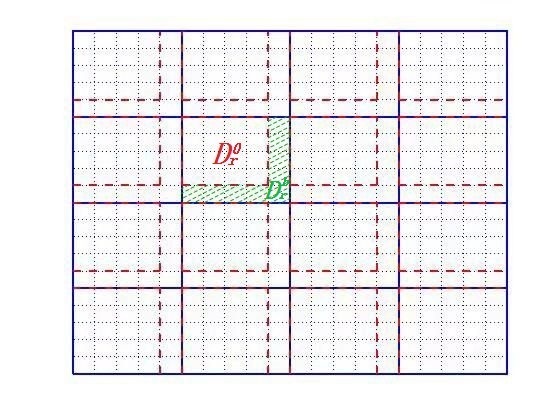}
\caption{The left graph: each small square with dotted lines denotes an element $\Omega_k$ and each square with solid lines denotes a subdomain $D_r$.
The right graph: each rectangle with red dotted lines denotes a subdomain $D_r^0$, each $L$-shape or revere $L$-shape domain denotes a subdomain $D_r^b$ (the green shade domain), where $D_r^b\cup D_r^0=D_r$.}
\label{dd2}
\end{figure}

%\subsection{A space decomposition}
As usual we coarsen the partition as follows: let $\Omega$ be decomposed into a union of $D_1,D_2,\ldots,D_{n_0}$ such that $D_r$ is just a union of several elements $\Omega_k \in {\mathcal T}_h$
%(see Figure \ref{dd2})
and satisfies (refer to the left graph of Figure \ref{dd2})
$$ \overline{\Omega} = \bigcup\limits_{r=1}^{n_0}\overline{D}_r, \quad D_r\bigcap D_l = \emptyset \quad \text{for } r\ne l. $$
Let $d$ denote the size of the subdomains $D_1,D_2,\cdots,D_{n_0}$, and let
$\mathcal{T}_d$ denote the partition comprised of the subdomains $\{D_r\}_{r=1}^{n_0}$.% and $V_h^p(\mathcal{T}_d)$ denote the space consisting of polynomial functions associated with the coarse triangulation $\mathcal{T}_d$.

For the construction of a substructuring preconditioner, we need to define a suitable ``interface" $\Gamma$ such that the degrees of freedoms in all the subdomain interiors (i.e., $\Omega\backslash \Gamma$)
can be eliminated independently for different subdomains.
We first explain that, for the current situation,  an interface $\Gamma$ cannot be defined in the standard manner, where $\Gamma$ is just a union of all the intersections of two neighboring subdomains.
To this end, we want to investigate basis functions associated with two neighboring subdomains $D_r$ and $D_l$, which have the non-empty common part $\partial D_r\cap\partial D_l$.
Let $e$ and $e'$ be two fine edges that satisfy $e\in \bar{D}_r\backslash(\partial D_r\cap\partial D_l)$ and $e'\in \bar{D}_l\backslash(\partial D_r\cap\partial D_l)$,
and let $\mu_e$ and $\mu_{e'}$ denote two basis functions on $e$ and $e'$ respectively. It can be checked that, if $e$ and $e'$ are close to $\partial D_r\cap\partial D_l$,
then $\mu_e$ and $\mu_{e'}$ still have coupling, i.e., $s_h(\mu_e,\mu_{e'})\not=0$. This means that, if the interface is defined
in the standard manner, namely, is defined as the union of all $\partial D_r\cap\partial D_l$, the degrees of freedom in subdomain interiors cannot be eliminated independently.
According to this observation, in the current situation an interface should be defined as a union of some elements instead of a union of some edges.

For each $D_r$, let $D_r^b\subset D_r$ be a union of the elements that touch the right and the lower boundary of $\partial D_r\backslash\partial\Omega$ (refer to the right graph in Figure \ref{dd2}).
We define an interface as
$$ \Gamma=\bigcup\limits_{r=1}^{n_0}D_r^b. $$
Of course, the definition of such an interface is not unique (see \cite{HuZ2016} and \cite{Shu2018} for similar definitions of interfaces), for example, an interface $\Gamma$ can be
defined as a union of all the elements that touch the standard interface $\cup_{k\not=j}(\partial D_k\cap\partial D_j)$.

In the following we describe various subspaces of $W_h^q(\gamma)$ and the corresponding solvers, which are needed in the construction of the desired preconditioner.

At first we define a subspace associated with each $D_r$. Set $D_r^0 = D_r\backslash D_r^b$ (see the right graph in Figure \ref{dd2}), and define the subspace for each subdomain $D^0_r$
$$ W_h^q(D_r^0) = \{  \mu\in W_h^q(\gamma):supp\,\mu\subset D_r^0 \}, ~ r=1,2,\dots,n_0.$$
The local solver on the local space $ W_h^q(D_r^0)$ is defined in the standard manner. Let $S_r^0: W_h^q(D_r^0)\to W_h^q(D_r^0)$ be the restriction of $S$ on $W_h^q(D_r^0)$
$$ \langle S_r^0\varphi, \psi\rangle=\langle S\varphi,\psi\rangle,\quad \varphi,\psi\in W_h^q(D_r^0). $$
%\begin{center}
%\end{center}

\begin{figure}[!ht]
\centering
\includegraphics[width=5.3cm,height=5cm]{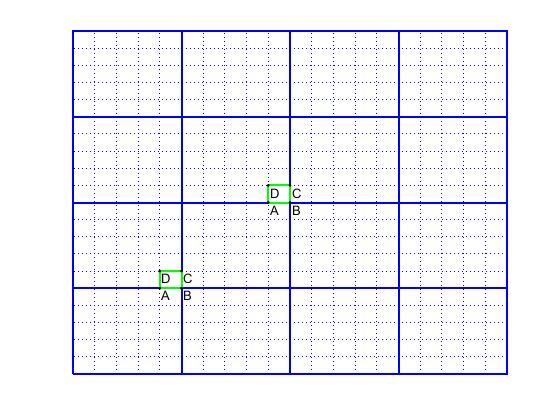}\quad\quad\includegraphics[width=5.3cm,height=5cm]{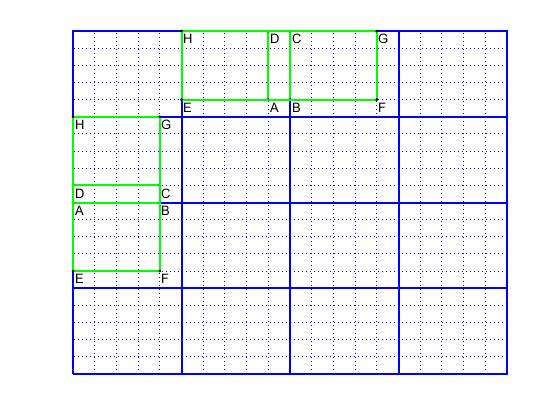}
\caption{The left graph: the rectangle ABCD denotes a subdomain $D_V$. The right graph: the rectangle ABCD denotes a subdomain $D_{rl}$ and the rectangle EFGH denotes a subdomain $\tilde{D}_{rl}$.}
\label{dd3}
\end{figure}

For the definition of solvers associated with the interface, we need to give a decomposition of the interface $\Gamma$.
Let $\mathcal{V}_d$ denote the set of all the nodes corresponding to the coarse partition $\mathcal{T}_d$. For a coarse node $V\in\mathcal{V}_d$, let $D_{V}$ denote the
top left corner element that touch the vertex $V$ (see the left graph in Figure \ref{dd3}).

Let $D_{rl}$ denote the union of the elements that touch the intersection $\partial D_r \cap \partial D_l$ from the left side (or the upper side) but do not touch the lower (or the right) endpoints of
$\partial D_r\cap\partial D_l$ (see the right graph in Figure \ref{dd3}).

It is easy to see that the interface can be decomposed into
$$ \Gamma=\big(\bigcup\limits_{rl} D_{rl}\big)\bigcup\big(\bigcup\limits_{V\in\mathcal{V}_d}D_V\big). $$

Next we define local interface spaces. Set
\begin{equation*}
\tilde{D}_{rl}=D_{rl}\cup D_r^0\cup D_l^0.
\end{equation*}
and define the discrete $s_h(\cdot,\cdot)$-harmonic extension spaces
$$ W_H^q(\tilde{D}_{rl}) = \{ \mu\in W_h^q(\gamma):supp\,\mu\subset \tilde{D}_{rl}; s_h(\mu,w)=0,\forall w\in W_h^q(D_r^0)\cup W_h^q(D_l^0)\}. $$

Notice that the basis functions of these local spaces are not given explicitly, so the variational problems defined on these spaces cannot be solved in the direct manner.
In order to overcome this difficulty, instead of computing such basis functions, as usual (see, for example, \cite{bddc}) we transform the corresponding local interface problem into a residual equation
, which is defined on the natural restriction space of the global space $W_h^q(\gamma)$ on the subdomain $\tilde{D}_{rl}$ (such residual equation
will be described exactly in Step 2 of Algorithm 3.1).
However, solution of the residual equation is expensive since the restriction space contains many more basis functions than each local space $W_h^q(D_r^0)$, which is defined on a smaller subdomain $D_r^0$ than $\tilde{D}_{rl}$.

\begin{figure}[!ht]
\setlength{\abovecaptionskip}{-8pt}
\setlength{\belowcaptionskip}{-8pt}
\centering
\includegraphics[width=6cm,height=5.7cm]{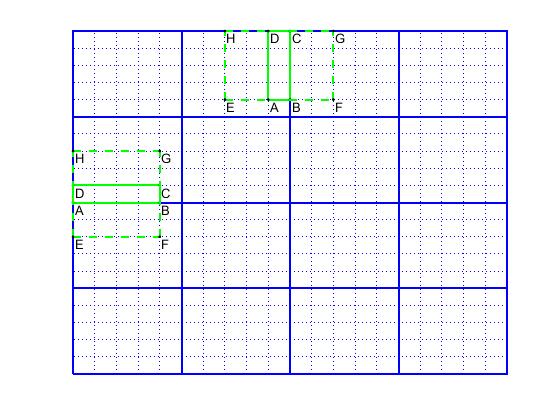}
\caption{The rectangle ABCD denotes a subdomain $D_{rl}$ and the rectangle EFGH denotes a subdomain $\tilde{D}_{rl}^{half}$}
\label{dd5}
\end{figure}

In order to decrease the cost of calculation, we choose to reduce the sizes of the subdomains $\tilde{D}_{rl}$ and define discrete $s_h(\cdot,\cdot)$-harmonic on the reduced subdomains. We reduce $\tilde{D}_{rl}$ to $\tilde{D}_{rl}^{half}$ such that the resulting subdomains  have almost the same size $d$ with $D_r$ (see Figure \ref{dd5}).

Define the local spaces
$$ W_h^q(\tilde{D}_{rl}^{half})=\{ \mu\in W_h^q(\gamma):supp\,\mu\subset \tilde{D}_{rl}^{half}\}.$$
For $\mu\in W_H^q(\tilde{D}_{rl})$, define $\mu_{rl}^{half}\in W_h^q(\tilde{D}_{rl}^{half})$ such that $\mu_{rl}^{half}|_{D_{rl}}=\mu|_{D_{rl}}$ and
$\mu_{rl}^{half}$ is discrete $s_h(\cdot,\cdot)$-harmonic in the complement domain $\tilde{D}_{rl}^{half}\backslash D_{rl}$.

Define the discrete operator $K_{rl}^0: W_H^q(\tilde{D}_{rl})\to W_H^q(\tilde{D}_{rl})$ by
$$ \langle K_{rl}^0 \mu, w\rangle=s_h(\mu_{rl}^{half},w_{rl}^{half}), ~\mu\in W_H^q(\tilde{D}_{rl}), ~\forall w\in W_H^q(\tilde{D}_{rl}).$$

Notice that the action of $(K_{rl}^0)^{-1}$ is implemented by solving a residual equation defined on the ``half'' space $W_h^q(\tilde{D}_{rl}^{half})$ (see Algorithm 3.1),
so $K_{rl}^0$ can be regarded as an ``inexact'' local interface solver based on the ``compressed'' harmonic extension $\mu_{rl}^{half}$ (refer to \cite{HuZ2016}).
It is easy to see that the dimension of $W_h^q(\tilde{D}_{rl}^{half})$ is about half of the dimension of $W^q_h(\tilde{D}_{rl})$ and almost equals the dimension of
$W_h^q(D^0_r)$. Then almost the same cost is needed for the solution of each subproblem in Step 2 and Step 1 (and Step 3) of Algorithm 3.1,
which make the loading balance be guaranteed in parallel calculation (in applications, we choose $d\approx\sqrt{h}$).

Finally we construct a coarse space $W_d^q(\gamma)$ by some local energy minimizations.

For a coarse node $V\in{\mathcal V}_d$, let $\phi^{(m)}_{V}$ be a basis function in the subspace
$$ W_h^q(D_V)= \{\mu|_{D_V}: \mu\in W_h^q(\gamma) \}. $$
Since the function $\phi^{(m)}_{V}$ is well defined only on the fine edges of $D_V$, we need to extend $\phi^{(m)}_{V}$
in a suitable manner such that $\phi^{(m)}_{V}$ has definitions on all the fine edges of ${\mathcal T}_h$. The desired coarse space will be spanned by the extensions of all $\phi^{(m)}_{V}$.

Let $\tilde{\phi}^{(m)}_{V}$ be the initial extension of $\phi^{(m)}_{V}$ such that $\tilde{\phi}^{(m)}_{V}$ is $s_h(\cdot,\cdot)$-harmonic on each subspace $W_h^q(D_r^0)$ and vanishes on all the
fine edges in $\Gamma\backslash D_V$. In order to define further extension of $\tilde{\phi}^{(m)}_{V}$, let $\Gamma_{V}$ denote a union of the coarse edges that touch the vertex $V$. For each $\Gamma_{rl}\in\Gamma_V$,
let $\Phi^{(m)}_{V,rl}\in W_h^q(\tilde{D}^{half}_{rl})$ be the solution of the
%extension such that $\tilde{\phi}^{(m)}_{V}+\Phi^{(m)}_{V,rl}$ solve the
minimization problem
\begin{equation}
\min\limits_{\Psi\in W_h^q(\tilde{D}^{half}_{rl})} \{s_h^{(r)}(\tilde{\phi}^{(m)}_{V}+\Psi, \tilde{\phi}^{(m)}_{V}+\Psi)+s^{(l)}_h(\tilde{\phi}^{(m)}_{V}+\Psi, \tilde{\phi}^{(m)}_{V}+\Psi)\},\label{3.min}
\end{equation}
where $s^{(r)}_h(\cdot,\cdot)$ denotes the restriction of $s_h(\cdot,\cdot)$ on the fine edges on $D_r$.
Then $\Phi^{(m)}_{V,rl}\in W_h^q(\tilde{D}^{half}_{rl})$ can be obtained by solving the local equation
\begin{equation}
\sum\limits_{k=r,l}s^{(k)}_h(\Phi^{(m)}_{V,rl}, v)
=-\sum\limits_{k=r,l}s^{(k)}_h(\tilde{\phi}^{(m)}_{V}, v),\quad\forall v\in W_h^q(\tilde{D}^{half}_{rl}).
\label{3.local}
\end{equation}
Define
\begin{equation}
\Phi^{(m)}_V = \tilde{\phi}^{(m)}_{V} + \sum_{\Gamma_{rl}\in\Gamma_V} R^t_{rl}\Phi^{(m)}_{V,rl},\label{3.coarse basis}
\end{equation}
where $R^t_{rl}$ denotes the zero extension operators from $W_h^q(\tilde{D}_{rl})$ into $W_h^q(\gamma)$.
The coarse space $W_d^q(\gamma)$ is spanned by all the basis functions $\Phi^{(m)}_V$, namely,
$$ W_d^q(\gamma) = \text{span}\{ \Phi^{(m)}_V \}. $$

%Let all basic functions of $W_q(D_V^r)$ be expanded to subspace $W_q(\bar{D}_r)$ by the harmonic expansion. Suppose there are m basic functions in space $W_q(D_V^r)$, $\{\phi_k^r, ~k = 1,\cdots,m\}$. Let $\Phi_k^r$ be the harmonic expansion function of $\phi_k^r$ to $W_q(\bar{D}_r)$.

%Firstly, for each $r$, Let all basic functions of $W_q(D_V^r)$ be expanded to subspace $W_q(\bar{D}_r)$ by the harmonic expansion. Suppose there are m basic functions in space $W_q(D_V^r)$, $\{\phi_k^r, ~k = 1,\cdots,m\}$. Let $\Phi_k^r$ be the harmonic expansion function of $\phi_k^r$ to $W_q(\bar{D}_r)$.

%We notice that every subdomain $D_V$ belongs to four sundomains $\bar{D}_r$, which can be written as $\bar{D}_1$, $\bar{D}_2$, $\bar{D}_3$ and $\bar{D}_4$. A basic function $\phi_k(supp\,\phi_k \subset D_V)$ can be expanded to $\Phi_k^1$, $\Phi_k^2$, $\Phi_k^3$ and $\Phi_k^4$ by harmonic expansion and then make the weighted sum and zero extension to the whole discrete space, which implies
%$$ \Phi_k = \sum_{r=1}^4 (\bar{R}_r)^t {\mathcal Q}_r \Phi_k^r, $$

Let the coarse solver $S_d: W_d^q(\gamma)\to W_d^q(\gamma)$ be the discrete operator which is the restriction of $S$ on $W_d^q(\gamma)$ as usual.

%Let $R_d$ denote the restriction operators from space $W_h^q(\gamma)$ to subspaces $W_d^q(\gamma)$. Let $A_d$ denote the stiffness matrix induced from the sesqulinear form $a_h(\cdot,\cdot)$ on the subspace $W_d^q(\gamma)$.
Now we can define the preconditioner $K: W_h^q(\gamma)\to W_h^q(\gamma)$ as
$$ K^{-1}=\sum_r (S_r^0)^{-1}Q_r + \sum_{\Gamma_{rl}} (K_{rl}^0)^{-1}Q_{rl} + S_d^{-1}Q_d, $$
where $Q_r, Q_{rl}$ and $Q_d$ denote the $L^2$ projectors into $W_h^q(D_r^0), W_H^q(\tilde{D}_{rl})$ and $W_d^q(\gamma)$, respectively.

The action of the preconditioner $K^{-1}$ can be described by the following algorithm.

{\bf Algorithm 3.1}. For $\xi\in W_h^q(\gamma)$, the solution $\lambda_{\xi}= K^{-1}\xi\in W_h^q(\gamma)$ can be obtained as follows:

Step 1. Computing $\lambda_r^0\in W_h^q(D_r^0)$ in parallel by
%$$ \lambda_1 = \sum_{k=1}^{n_0} (R_k^0)^t(A_k^0)^{-1}R_k^0g. $$
$$  s_h^{(r)}(\lambda_r^0, \mu_h) =\langle\xi, \mu_h\rangle, \quad\forall\mu_h\in W_h^q(D_r^0), ~r=1,2,\dots,n_0. $$

Step 2. Computing $\lambda_{rl}\in W_h^q(\tilde{D}_{rl}^{half})$  in parallel by
%$$ \lambda_2 = \sum_{\Gamma_{kj}} (R_{kj})^t(A_{kj})^{-1}R_{kj}(g-A\lambda_1). $$
$$  \sum\limits_{k=r,l} s^{(k)}_h(\lambda_{rl}, \mu_h) = \langle\xi, \mu_h\rangle-\sum\limits_{k=r,l}s^{(k)}_h(\lambda_r^0, \mu_h), \quad\forall\mu_h\in W_h^q(\tilde{D}_{rl}^{half}).$$

Step 3. Computing $\lambda_d\in W_d^q(\gamma)$ by
%$$ \lambda_3 = (R_d)^t(A_d)^{-1}R_d(g-A\lambda_1). $$
$$  s_h(\lambda_d, \mu_h) = \langle\xi, \mu_h\rangle-\sum_{r} s^{(r)}_h(\lambda_r^0, \mu_h), \quad\forall\mu_h\in W_d^q(\gamma).$$

Step 4. Set $\phi = \sum\lambda_{rl} + \lambda_d$, and compute
%$\phi = \lambda_2+\lambda_3$
harmonic extensions $\lambda^H_r\in W_h^q(D_r)$ for all $r$ in parallel, such that $\lambda^H_r = \phi$ on $D^b_r$ and satisfies
$$ s^{(r)}_h(\lambda^H_r, ~\mu_h)=0,\quad \forall \mu_h\in W_h^q(D_r^0), ~r=1,2,\cdots,n_0.$$

Step 5. Computing
%$$\lambda_g=\lambda_1+\sum_{r}\lambda^H_r.$$
$$\lambda_{\xi} = \sum_{r} \lambda_r^0 + \sum_{r}\lambda^H_r.$$
\begin{remark}
The minimization problem (\ref{3.min}) is different from that in the BDDC method. In the BDDC method, each minimization problem which determines coarse basis functions
is defined on one subdomain, so the solutions of the two minimization problems associated with two neighboring subdomains have different values on their common interface.
In order to define coarse basis functions, in the BDDC method one has to compute some average of the values of the two solutions on the common interface. Since the minimization problem
(\ref{3.min}) is defined on the subdomain $\tilde{D}^{half}_{rl}$, the solution of this minimization problem has a unique value on the interface $D_{rl}$ and the coarse
basis functions can be directly obtained by (\ref{3.coarse basis}). We found that, if minimization problems are defined as in the BDDC method, then the resulting preconditioner
is unstable.
\end{remark}

\begin{remark} Since the stiffness matrix of $S^r_h$ has almost the same structure as the stiffness matrix ${\mathcal S}$ of the global system, the condition number of the stiffness matrix of $S^r_h$
cannot be significantly decreased comparing the original stiffness matrix ${\mathcal S}$. However, such a local stiffness matrix has much lower order than
${\mathcal S}$, so each subproblem in Step 1 of Algorithm 3.1 can be solved in a direct manner (using LU decomposition), which is not sensitive to the condition number of this local stiffness matrix,
where the global Step 1 is implemented in parallel.
Notice that the variational problem (\ref{3.local}) and the variational problem in Step 2 of Algorithm 3.1 correspond to the same stiffness matrix
(with different right hands only). Thus the computation for the coarse basis functions by solving every subproblem (\ref{3.local}) in parallel only increases a
little cost by using LU decomposition made in Step 2 for each local stiffness matrix (when Step 2 is implemented in the direct method).
\end{remark}

\begin{remark} When $\Omega$ is a general domain than a rectangle, we can first define a domain decomposition such that every subdomain $D_r$ is a polygon, and then define a triangle partition
on each subdomain $D_r$, all of which constitute a partition ${\mathcal T}_h$ of $\Omega$. In this situation, the ``interface" $\Gamma$ and the reduced subdomain $\tilde{D}_{rl}^{half}$ can be
defined in a similar manner, but their shapes may be more complicated.
\end{remark}

\section{Main results} \label{theo}
Throughout this paper, $C$ denotes a generic positive constant that may have different values in different occurrences, where $C$ is always independent of $\omega, h,p$ and $q$
but may depend on the shape of $\Omega$ and the maximal value and minimal value of $c({\bf x})$ on $\Omega$. Before presenting the main results, we give several assumptions. \\
{\bf Assumption 1}. The domain $\Omega$ is a strictly star-shaped; the function $c({\bf x})$ belongs to $W^{1,\infty}(\Omega)$.

The first condition in the above assumption appeared in many existing works to build error estimates with little wave number pollution (see, for example, \cite{DuWu2015} and \cite{fem});
the second condition was used in \cite{Brown2016} to build stability result of analytic solution. \\
{\bf Assumption 2}. The mesh size $h$ satisfies the condition: $\omega h\leq C_0$ with a possibly small constant $C_0$ independent of
$\omega$, $h,~p$ and $q$ (but  may depend on the shape of $\Omega$ and the maximal value and minimal value of the function $c({\bf x})$).

The above assumption is weaker than that required in analysis of the HDG-type methods. The following assumption has no restriction to the proposed method.\\
{\bf Assumption 3}. The parameter $\rho$ in the variational formula is not large: $\rho\leq C_0\min\{1,\omega^2h\}$ for a possibly small constant $C_0$ independent of
$\omega$, $h,~p$ and $q$.

From the viewpoint of algorithm, all the discretization methods based on polynomial basis functions are practical for the case with variable wave numbers (in inhomogeneous media).
However, almost existing error estimates with little wave number pollution were established only for the case with constant wave numbers (see, for example, \cite{DuWu2015} and \cite{fem}).
The main reason is that one does not know whether the result on ``stable decomposition of solution", which was built in Theorem 4.10 of \cite{fem} and plays a key role in the derivations of good error estimates,
still holds for the case with variable wave numbers. In this paper we try to investigate the possibility that the proposed method possesses error estimates with little wave number pollution
even for the case of variable wave numbers. In order to cover the case of variable wave number, we introduce an additional assumption.

For $\tilde{f}\in L^2(\Omega)$, consider a dual problem with Robin-type boundary condition
\begin{equation} \label{adjo1}
\left\{
\begin{aligned}
& -\Delta\phi - \kappa^2 \phi = \tilde{f}
& \text{in}\ \Omega,\\
& \frac{\partial \phi}{\partial n} - i\kappa \phi = 0
& \text{on}\ \partial\Omega.
\end{aligned}
\right.
\end{equation}
Let $\tilde{V}_h^p({\mathcal T}_h)\subset H^1(\Omega)$ denote the continuous piecewise $p$-order polynomial space associated with the partition ${\mathcal T}_h$.\\
{\bf Assumption 4}. The finite element solution $\phi_h\in \tilde{V}_h^p({\mathcal T}_h)$ of (\ref{adjo1}) possesses a weak convergence with respect to $p$ for large $p$
\begin{equation}
||\nabla(\phi-\phi_h)||_{0,\Omega} + \omega||(\phi-\phi_h)||_{0,\Omega}
\lesssim p^{-{1\over 2}}||\tilde{f}||_{0,\Omega}. \label{4.assumption4}
\end{equation}

This assumption can be met easily when $c({\bf x})$ is a constant. In fact, for this case the following stronger result has been built in \cite[Cor 5.10]{fem} under the
assumptions that $\Omega$ is a strictly star-shaped domain with an analytic boundary and the discretization parameters satisfy the mild conditions ${\omega h\over p}\leq C_0$ and $p\geq 1+c_0\log\omega$:
\begin{equation}
||\nabla(\phi-\phi_h)||_{0,\Omega} + \omega||(\phi-\phi_h)||_{0,\Omega}
\lesssim hp^{-1} ||\tilde{f}||_{0,\Omega}. \label{5.2new06}
\end{equation}
Therefore, when $c({\bf x})$ is a constant, Assumption 4 should be changed into: $\Omega$ has an analytic boundary; $p\geq 1+c_0\log\omega$ (since Assumption 2 implies ${\omega h\over p}\leq C_0$).

Whether the error estimate (\ref{5.2new06}) still holds for the case of variable $c({\bf x})$ seems an open problem, but the weak error estimate (\ref{4.assumption4}) should be valid even for a 
variable $c({\bf x})$ under the above assumptions. In this situation, Assumption 4 can be replaced by the conditions that $\Omega$ has an analytic boundary and $p\geq 1+c_0\log\omega$.

Now we list the main results, which will be proved in the next section.
Firstly, we give a result about local $inf-sup$ condition.

\begin{theorem} \label{infsup} Let $q\geq 1$ and $p\geq q+2$. For any $\mu\in W_h^q(\partial\Omega_k \backslash\partial\Omega)$,
there exits a non-zero function $v\in V_h^p(\partial\Omega_k)$ such that
\begin{equation} \label{is}
\langle\mu, v\rangle_{\partial\Omega_k \backslash\partial\Omega}
\geq C q^{-\frac{1}{2}} ||\mu||_{0,\partial\Omega_k \backslash\partial\Omega} ||v||_{0,\partial\Omega_k\backslash\partial\Omega},
\end{equation}
where $C$ is a constant independent of $\omega,~h~,p$ and $q$.
\end{theorem}

Next we give a result on the coerciveness of the sesquilinear form $s_h(\cdot,\cdot)$, which implies that the discrete problem (\ref{discre}) is well posed.

\begin{theorem} \label{well} Let Assumption 1-Assumption 4 be satisfied. Suppose $q\geq 1$ and  $p\geq q+2$. Then, for any $\mu_h\in W_h^q(\gamma)$, we have
\begin{equation}
 s_h(\mu_h,\mu_h)
\geq C \omega^{-2}h^{2}p^{-1}q^{-1} \sum_{\gamma_{kj}} \|\mu_h||^2_{0,\gamma_{kj}}, \label{4.new1}
\end{equation}
where $C$ is a constant independent of $\omega,~h~,p$ and $q$.
\end{theorem}

Finally, we give error estimates of the approximation $u_h$.
Define the subspace $W^{r-\frac{1}{2}}(\gamma) = \prod_{k\not=j} H^{r-\frac{1}{2}}(\gamma_{kj})$, which is equipped with the norm
$$||\mu||_{r-{1\over 2},\gamma}=\big(\sum_{\gamma_{kj}} ||\mu||^2_{r-\frac{1}{2},\gamma_{kj}}\big)^{{1\over 2}},\quad\forall~\mu\in W^{r-\frac{1}{2}}(\gamma).$$
For ease of notation, we also define the semi-norms ($r\geq 1$)
$$ |\mu|_{r-{1\over 2},\gamma}=\big(\sum_{\gamma_{kj}} |\mu|^2_{r-\frac{1}{2},\gamma_{kj}}\big)^{{1\over 2}},\quad~\mu\in W^{r-\frac{1}{2}}(\gamma) $$
and
$$ |v|_{r,\Omega}=\big(\sum\limits_{k=1}^N |v|^2_{r,\Omega_k}\big)^\frac{1}{2}, \quad~v\in \prod_{k=1}^NH^r(\Omega_k). $$
Set
$$ H^{r+1}({\mathcal T}_h)=\{v\in H^2(\Omega):~~v|_{\Omega_k}\in H^{r+1}(\Omega_k),~{\partial v\over \partial{\bf n}}|_{\gamma_{kj}}\in H^{r-{1\over 2}}(\gamma_{kj})\}. $$
%\mbox{for~each}~k, \mbox{for~each}~\gamma_{kj}
\begin{theorem} \label{fina} Suppose that $q\geq 1$ and $p\geq q+2$. Let Assumption 1-Assumption 4 be satisfied.
Assume that the analytical solution $u$ of the Helmholtz problem (\ref{eq}) belongs to $H^{r+1}({\mathcal T}_h)$ with $1\leq r\leq q$ ($r\in \mathbb{N}$).
Then the approximate solution $u_h$ defined in Subsection \ref{uh} satisfies
\begin{equation}
|u-u_h|_{1,\Omega}\leq C h^{r-1} \big(p^{-r}|u|_{r+1,\Omega} +q^{-r}|\lambda|_{r-{1\over 2},\gamma}\big)\label{4.new2}
\end{equation}
and
\begin{equation}
||u-u_h||_{0,\Omega}
\leq C \omega^{-1}h^{r-1} \big(p^{-r}|u|_{r+1,\Omega}+q^{-r}|\lambda|_{r-{1\over 2},\gamma} \big), \label{4.new3}
\end{equation}
where $C$ is a constant independent of $\omega,~h~,p$ and $q$.
\end{theorem}

\begin{remark} Comparing Theorem \ref{fina} with Theorem 3.15 in \cite{Hiptmair2011} (and  Theorem 3.4 in \cite{Monk1999}), we can see that the proposed discretization method possesses
almost the same $L^2$ convergence order as the plane wave methods (for the case of constant wave number), which have fast convergence and small ``wave number pollution".
As pointed out in Section 1, the standard plane wave methods are not practical for the case with variable wave numbers, but there is not this problem for the proposed method (see Remark 2.1).
\end{remark}

\section{Proof of the main results}
In this section, we give detailed proofs of the theorems stated in Section \ref{theo}. Since the proposed approximate solution $(u_h,\lambda_h)$ do not satisfy a mixed variational problem
(as in the Lagrange multiplier methods) or a hybridizable variational problem (as in the hybridizable discontinuous Galerkin methods), so the results cannot be proved by the techniques
developed in existing works.
As we shall see, the proofs are very technical, so we divide this section into three subsections, in which we will build many auxiliary results.

For ease of notation,
we use the shorthand notation $x \lesssim y$ and $y \gtrsim x$ for the inequality $x \leq Cy$ and $y \geq Cx$, where $C$ is a constant independent of $\omega$, $h$, $p$ and $q$ but
may depend on the shape of $\Omega$ and the maximal value and minimal value of $c({\bf x})$ on $\Omega$.
Throughout this this section, we use $p$ and $q$ to denote two positive integers.

We first verify the local $inf-sup$ condition given in Theorem 4.1 by using Jacobi polynomials.

\subsection{Analysis on the local $inf-sup$ condition} \label{51}
The main difficulty for the proof of Theorem 4.1 is the fact that the functions in $W_h^q(\partial\Omega_k\backslash\partial\Omega)$ are defined independently for different edges of $\partial\Omega_k$
and may be discontinuous at the vertices of $\partial\Omega_k$ but the functions in $V_h^p(\partial\Omega_k)$ are defined globally on $\partial\Omega_k$ and must be continuous at
the vertices of $\partial\Omega_k$. Because of this, we have to
%separately consider every edge of $\partial\Omega_k$ and
split the $q$-order polynomial space into a sum of two polynomial subspaces,
one of which consists of all the $q$-order polynomials vanishing at the vertices of $\partial\Omega_k$, so that the construction of a function $v$ satisfying (\ref{is}) for $\mu\in W_h^q(\partial\Omega_k\backslash\partial\Omega)$ becomes easier by using this splitting and Jacobi polynomial basis functions of this subspace (we can require that such function $v$ vanishes at all the vertices of $\partial\Omega_k$).

Set $J=[0,1]$ and let ${\mathcal P}_q$ stands for the space of all polynomials on $J$ with orders $\leq q$.
Firstly, we give a space decomposition of ${\mathcal P}_q$ on $J$
\begin{equation} \label{com}
{\mathcal P}_q = {\mathcal P}_1^*+{\mathcal P}'_{q}, \text{ with } {\mathcal P}_1^*\perp{\mathcal P}'_{q}.
\end{equation}
The specific definition of ${\mathcal P}_1^*$ and ${\mathcal P}'_{q}$ will be given next.

Let ${\mathcal P}_1$ and ${\mathcal P}'_{q}$ denote the linear part and high-order part of ${\mathcal P}_{q}$, respectively. Then the two basis functions of ${\mathcal P}_1$ are $\phi_1=x$ and $\phi_2=1-x$.
If $q=1$, then ${\mathcal P}'_{q}=\emptyset$ and set $P_1^*=\{\phi_1, \phi_2\}$. For a unified description below, we define $\phi_1^*=\phi_1$ and $\phi_2^*=\phi_2$ when $q=1$, and write $P_1^*=\{\phi^*_1, \phi_2^*\}$.

In the following we assume that $q\geq 2$. Let $\{\psi_k\}_{k=1}^{q-1}$ denote the basis functions of the subspace ${\mathcal P}'_{q}$. Define
\begin{equation}
\phi_1^* = \phi_1-\sum_{k=1}^{q-1} \alpha_k\psi_k\quad\mbox{and}\quad
\phi_2^* = \phi_2-\sum_{k=1}^{q-1} \beta_k\psi_k. \label{5.newequal}
\end{equation}
Here the numbers $\{\alpha_k\}$ and $\{\beta_k\}$ are determined by $\langle\phi_1^*,\psi_k\rangle_J=0$ and $\langle\phi_2^*,\psi_k\rangle_J=0$.
Apparently we can get
\begin{equation*}
(\alpha_1, \alpha_2, \dots, \alpha_{q-1})^t = A^{-1}b_1
\text{ and }
(\beta_1, \beta_2, \dots, \beta_{q-1})^t = A^{-1}b_2,
\end{equation*}
where $A \!=\! (\langle \psi_k,\psi_j\rangle_J)_{(q-1)*(q-1)}$ and $b_1 \!=\! (\langle \phi_1,\psi_k\rangle_J)_{(q-1)*1},b_2 \!=\! (\langle \phi_2,\psi_k\rangle_J)_{(q-1)*1}$,
which means $\phi_1^*, \phi_2^*$ are uniquely determined.
Let ${\mathcal P}_1^*=\text{span} \{\phi_1^*, \phi_2^*\}$, which satisfies the space decomposition (\ref{com}).

Next we give a set of orthogonal basis functions of ${\mathcal P}'_{q}=\text{span}\{\psi_1,\psi_2,\dots,\psi_{q-1}\}$ ($q\geq 2$). In order to explicitly
write the orthogonal basis functions and conveniently compute the involved integrations, we use a set of Jacobi polynomials $\{G_k\}$ (see \cite{orth}).
For convenience, we let the coefficient of the first Jacobi polynomial be $1$. Then, for $p\geq q+2$, the Jacobi polynomials are defined as
\begin{equation}
G_k = (-1)^{k-1}\frac{(k+3)!}{(2k+2)!} x^{-2}(1-x)^{-2}\frac{d^{k-1}}{dx^{k-1}}(x^{k+1}(1-x)^{k+1}),\quad k=1,\cdots,p.
\end{equation}
It is known that
\begin{equation*}
\int_0^1 x^2(1-x)^2 G_k G_j dx=
\begin{cases}
\quad 0 &  k \ne j, \\
\frac{(k-1)!(k+1)!^2(k+3)!}{(2k+2)!(2k+3)!}  &  k = j.
\end{cases}
\end{equation*}
We also have the recursion relations
\begin{equation*}
\left\{
\begin{aligned}
& G_1 = 1, \ G_2 = x-\frac{1}{2}, \\
& G_k = (x-\frac{1}{2})G_{k-1}-\frac{(k-2)(k+2)}{4(2k-1)(2k+1)}G_{k-2}, \quad 3\leq k\leq p.
\end{aligned}
\right.
\end{equation*}

Define
\begin{equation}
\psi_k = \frac{(2k+2)!}{(k-1)!(k+1)!}x(1-x)G_k,\quad k=1,\cdots,q,\cdots,p. \label{5.new}
\end{equation}
It is clear that $\psi_k(0)=\psi_k(1)=0$ and
\begin{equation} \label{orth}
\int_0^1 \psi_k \psi_j dx=
\begin{cases}
\quad\quad 0 &  k \ne j, \\
\frac{k(k+1)(k+2)(k+3)}{2k+3}  &  k = j.
\end{cases}
\end{equation}
Furthermore $\{\psi_k\}_{k=1}^{p-1}$ satisfy the recursion relations
\begin{equation} \label{rec}
\left\{
\begin{aligned}
& \psi_1 = 12x(1-x), \ \psi_2 = 120x(1-x)(x-\frac{1}{2}), \\
& \psi_k = \frac{2(2k+1)}{k-1}(x-\frac{1}{2})\psi_{k-1}-\frac{k+2}{k-1}\psi_{k-2}, \quad 3\leq k\leq p.
\end{aligned}
\right.
\end{equation}
The functions $\{\psi_k\}_{k=1}^{q-1}$ constitute a set of orthogonal bases of ${\mathcal P}'_{q}$.

\begin{lemma} \label{ditui} Let $q\geq 2$.
For $\phi_1=x$, $\phi_2=1-x$ and $\{\psi_k\}_{k=1}^{q-1}$ defined by (\ref{5.new}), we have
\begin{equation*}
\langle \phi_1,\psi_k\rangle_J = 1
\text{ and }
\langle \phi_2,\psi_k\rangle_J = (-1)^{k-1}, ~k = 1,2,\dots,q-1.
\end{equation*}
\end{lemma}

\paragraph*{Proof.}
Using mathematical induction and the recursion relations (\ref{rec}), we can easily prove
\begin{equation*}
\langle \phi_1,\psi_k\rangle_J
=\int_0^1 x\psi_k dx = 1
\end{equation*}
and
\begin{equation*}
\langle \phi_2,\psi_k\rangle_J
=\int_0^1(1-x)\psi_k dx = (-1)^{k-1}.
\end{equation*}
$\hfill \Box$

It is easy to see that the following two equalities hold
%\begin{lemma}
for any positive integer $m$
%we have
\begin{equation} \label{gu1}
\sum_{k=1}^m \frac{2k+3}{k(k+1)(k+2)(k+3)}
=\frac{1}{3} - \frac{1}{(m+1)(m+3)}
\end{equation}
and
\begin{equation} \label{gu2}
\sum_{k=1}^m \frac{(-1)^{k-1}(2k+3)}{k(k+1)(k+2)(k+3)}
= \frac{1}{6} +\frac{(-1)^{m-1}}{(m+1)(m+2)(m+3)}.
\end{equation}
%\end{lemma}

%\paragraph*{Proof.} The two equalities are direct consequences of the following splitting
%$$
%\sum_{k=1}^m \frac{2k+3}{k(k+1)(k+2)(k+3)}
%=\sum_{k=1}^m \big(\frac{1}{k(k+2)}-\frac{1}{(k+1)(k+3)}\big)  $$
%and
%$$ \sum_{k=1}^m \frac{(-1)^{k-1}(2k+3)}{k(k+1)(k+2)(k+3)}
%=\sum_{k=1}^m  \big(\frac{(-1)^{k-1}}{k(k+1)(k+2)}+\frac{(-1)^{k-1}}{(k+1)(k+2)(k+3)}\big).
%$$
%$\hfill \Box$

\begin{lemma} \label{star} Let $\phi^*_1=x$ and $\phi_2^*=1-x$ for $q=1$. For $q\geq 2$, let $\{\phi_1^*, \phi_{2}^*\}$ be defined as (\ref{5.newequal}) with $\psi_k$
given by (\ref{5.new}). Then we have
\begin{equation} \label{phist}
\langle \phi_1^*,\phi_1^*\rangle_J
= \langle \phi_2^*,\phi_2^*\rangle_J
= \frac{1}{q(q+2)}, \quad
\langle \phi_1^*,\phi_2^*\rangle_J
=\frac{(-1)^{q-1}}{q(q+1)(q+2)}
\end{equation}
and
\begin{equation} \label{phistar}
||a_1\phi_1^*||^2_{0,J}+||a_2\phi_2^*||^2_{0,J}
\leq \frac{q+1}{q} ||a_1\phi_1^*+a_2\phi_2^*||^2_{0,J},
\quad\forall\, a_1,a_2 \in {\mathbb R}.
\end{equation}
\end{lemma}

\paragraph*{Proof.}

Using (\ref{5.newequal}), together with (\ref{orth}), (\ref{gu1}) and Lemma \ref{ditui}, we deduce that
\begin{equation*}
\begin{split}
\langle \phi_1^*,\phi_1^*\rangle_J
& %= \langle \phi_1,\phi_1^*\rangle_\Gamma
%& = \langle \phi_1,\phi_1\rangle_\Gamma -\big( \alpha_1\langle \phi_1,\psi_1\rangle_\Gamma+\alpha_2\langle \phi_1,\psi_2\rangle_\Gamma+\dots+\alpha_{p-1}\langle \phi_1,\psi_{p-1}\rangle_\Gamma \big) \\
= \langle \phi_1,\phi_1\rangle_J -b_1^tA^{-1}b_1
= \langle \phi_1,\phi_1\rangle_J -\sum_{k=1}^{q-1} \frac{\langle \phi_1,\psi_k\rangle^2_J}{\langle \psi_k,\psi_k\rangle_J} \\
& = \frac{1}{3} -\sum_{k=1}^{q-1} \frac{2k+3}{k(k+1)(k+2)(k+3)}
=\frac{1}{q(q+2)}.
\end{split}
\end{equation*}
Similarly, we have
\begin{equation*}
\langle \phi_2^*,\phi_2^*\rangle_J
=\frac{1}{q(q+2)}.
\end{equation*}

On the other hand, using Lemma \ref{ditui}, (\ref{orth}) and (\ref{gu2}), we get
\begin{equation*}
\begin{split}
\langle \phi_1^*,\phi_2^*\rangle_J
& %= \langle \phi_1,\phi_2^*\rangle_\Gamma
%& = \langle \phi_1,\phi_2\rangle_\Gamma -(\beta_1\langle \phi_1,\psi_1\rangle_\Gamma+\beta_2\langle \phi_1,\psi_2\rangle_\Gamma+\dots+\beta_{p-1}\langle \phi_1,\psi_{p-1}\rangle_\Gamma) \\
= \langle \phi_1,\phi_2\rangle_J-b_1^tA^{-1}b_2
= \langle \phi_1,\phi_2\rangle_J -\sum_{k=1}^{q-1} \frac{\langle \phi_1,\psi_k\rangle_J \,\langle \phi_2,\psi_k\rangle_J}{\langle \psi_k,\psi_k\rangle_J} \\
& = \frac{1}{6} -\sum_{k=1}^{q-1}\frac{(-1)^{k-1}(2k+3)}{k(k+1)(k+2)(k+3)}
= (-1)^{q-1}\frac{1}{q(q+1)(q+2)}.
\end{split}
\end{equation*}
This gives the second equality of (\ref{phist}). The equality (\ref{phistar}) can be easily obtained from (\ref{phist}).
$\hfill \Box$

%%%%%%%%%%%%%%%%%%%%%%%%%%%%%%%%%%%%%%%%%%%%%%%%%%%%%%%%%%%%%%%%%%%%%%%%%%%%%%
\paragraph*{{\bf Proof of Theorem \ref{infsup}}.}
For an element $\Omega_k$, let $n_k$ denote the number of the edges of $\Omega_k$ and write its boundary as $\partial \Omega_k = \bigcup _{j=1}^{n_k} J_j$, where $J_j$ is the $j$th edge of $\Omega_k$.
If $J_j\subset\partial\Omega_k\cap\partial\Omega$, we set $\mu|_{J_j} = 0$. Then we only need to prove: for any $\mu\in W_h^q(\partial\Omega_k)$, there exits a function $v\in V_h^p(\partial\Omega_k)$, such that
\begin{equation*}
\langle\mu, v\rangle_{\partial\Omega_k}
\geq C_{p,q} ||\mu||_{0,\partial\Omega_k} ||v||_{0,\partial\Omega_k},
\end{equation*}
where $C_{p,q}$ is a positive constant which may only depend on $p$ and $q$. Since $\Omega_k$ is regular, we can simply set $J_j=[0,1]$ by the scaling transformation.

By the space decomposition (\ref{com}), the function $\mu\in W_h^q(\partial\Omega_k)$ can be written as
\begin{equation}
\mu|_{J_j} =\sum_{k=1}^{q-1} \xi_{k} \psi_k+a_1\phi_1^*+a_2\phi_2^*,\quad a_1,a_2,\xi_{k}\in{\mathbb R}
\end{equation}
where $\{\phi_1^*,\phi_2^*\}$ are two basis functions of ${\mathcal P}_1^*$ and $\{\psi_k\}_{k=1}^{q-1}$ denote the orthogonal basis functions of ${\mathcal P}'_{q}$, see (\ref{5.new}).
Then we choose
\begin{equation}
v|_{J_j} =\sum_{k=1}^{q-1} \xi_{k} \psi_k+\sum_{k=q}^{p-1} \frac{\langle a_1\phi_1^*+a_2\phi_2^*,\psi_k \rangle_{J_j}}{\langle \psi_k,\psi_k \rangle_{J_j}} \psi_k,
\end{equation}
where $\{\psi_k\}_{k=1}^{p-1}$ are defined by (\ref{5.new}). It is clear that $v|_{J_j}(0)=v|_{J_j}(1)=0$. Then we have $v\in V_h^p(\partial\Omega_k)$.

Using the orthogonality condition (\ref{orth}), yields
\begin{equation*}
\langle\mu, v\rangle_{J_j}
= \sum_{k=1}^{q-1} \xi^2_{k} \langle\psi_k,\psi_k\rangle_{J_j}+\sum_{k=q}^{p-1} \frac{\langle a_1\phi_1^*+a_2\phi_2^*,\psi_k \rangle^2_{J_j}}{\langle \psi_k,\psi_k \rangle_{J_j}}
= ||v||^2_{0,J_j}.
\end{equation*}

It follows that
\begin{equation*}
\langle\mu, v\rangle_{\partial\Omega_k}
= \sum_{j=1}^{n_k} \langle\mu, v\rangle_{J_j}
= \sum_{j=1}^{n_k} ||v||^2_{0,J_j}
= ||v||^2_{0,\partial\Omega_k}.
\end{equation*}

Thus, we only need to prove: there exists $C_{p,q}$, such that
\begin{equation*}
||v||_{0,\partial\Omega_k}
\geq C_{p,q} ||\mu||_{0,\partial\Omega_k}
\quad\text{or}\quad
||v||_{0,J_j}
\geq C_{p,q} ||\mu||_{0,J_j}.
\end{equation*}

To do this, we use (\ref{orth}), (\ref{phist}) and Lemma \ref{ditui}, which gives
\begin{equation*}
\begin{split}
||v||^2_{0,J_j}
%& = \sum_{k=1}^{q-1} \xi^2_{k} \langle\psi_k,\psi_k\rangle_{\Gamma_j}+\sum_{k=q}^{p-1} \frac{\langle a_1\phi_1^*+a_2\phi_2^*,\psi_k \rangle^2_{\Gamma_j}}{\langle \psi_k,\psi_k \rangle_{\Gamma_j}} \\
& = \sum_{k=1}^{q-1} \xi^2_{k} \langle\psi_k,\psi_k\rangle_{J_j}+\sum_{k=q}^{p-1} \frac{(2k+3)(a_1+(-1)^{k-1}a_2)^2}{k(k+1)(k+2)(k+3)} \\
& = \sum_{k=1}^{q-1} \xi^2_{k} \langle\psi_k,\psi_k\rangle_{J_j}+\big(\frac{1}{q(q+2)}-\frac{1}{p(p+2)}\big) (a_1^2+a_2^2) \\
& +\big(\frac{(-1)^{q-1}}{q(q+1)(q+2)}-\frac{(-1)^{p-1}}{p(p+1)(p+2)}\big) 2a_1a_2.
\end{split}
\end{equation*}

Then, using (\ref{phist}) again, we have
\begin{equation*}
\begin{split}
||\mu||^2_{0,J_j}
& = \sum_{k=1}^{q-1} \xi^2_{k} \langle\psi_k,\psi_k\rangle_{J_j}+a_1^2\langle\phi_1^*,\phi_1^*\rangle_{J_j} +a_2^2\langle\phi_2^*,\phi_2^*\rangle_{J_j} +2a_1a_2\langle\phi_1^*,\phi_2^*\rangle_{J_j} \\
& = \sum_{k=1}^{q-1} \xi^2_{k} \langle\psi_k,\psi_k\rangle_{J_j}+\frac{1}{q(q+2)}(a_1^2+a_2^2) +\frac{(-1)^{q-1}}{q(q+1)(q+2)}2a_1a_2.
\end{split}
\end{equation*}

So we choose
\begin{equation*}
C^2_{p,q} =
\begin{cases}
1-\frac{(q+1)(q+2)}{(p+1)(p+2)}  & p+q = \text{even}, \\
1-\frac{(q+1)(q+2)}{p(p+1)}      & p+q = \text{odd},
\end{cases}
\end{equation*}
satisfying
\begin{equation*}
||u||^2_{0,J_j}\geq C^2_{p,q} ||\mu||^2_{0,J_j}.
\end{equation*}

Since $q\geq 1$ and $p\geq q+2$, we have $C^2_{p,q}\geq {2\over q+3}$ and so $||v||^2_{0,J_j}\gtrsim q^{-1} ||\mu||^2_{0,J_j}$, namely,
$$\langle\mu, v\rangle_{J_j}
\gtrsim q^{-\frac{1}{2}} ||\mu||_{0,J_j} ||v||_{0,J_j}.
$$
It concludes the proof of the local $inf-sup$ condition given by (\ref{is}).
$\hfill \Box$

\begin{remark} If $q\geq 2$ and $p=2q$, we have $C^2_{p,q}\geq {1\over 2}$, which implies that
$$
\langle\mu, v\rangle_{J_j}\gtrsim ||\mu||_{0,J_j} ||v||_{0,J_j}.
$$
Then, when $q\geq 2$ and $p=2q$, the inequality (\ref{is}) can be replaced by the optimal {\it inf-sup} condition
$$\langle\mu, v\rangle_{\partial\Omega_k \backslash\partial\Omega}
\geq C ||\mu||_{0,\partial\Omega_k \backslash\partial\Omega} ||v||_{0,\partial\Omega_k}. $$
\end{remark}

\subsection{Analysis on the coerciveness}

The proofs of Theorem 4.2 and Theorem 4.3 will depend on a {\it jump-controlled} stability estimate (which will be given by {\bf Proposition 5.1}). A technical tool for the derivation of
this stability estimate is a $Poincar\acute{e}$-type inequality given by Lemma \ref{dual}. In order to prove this $Poincar\acute{e}$-type inequality, we have to
develop a special technique: construct a globally continuous $p$-finite element function to ``approximate" a piecewise continuous $p$-finite element function and derive a corresponding
``approximate" result (Lemma \ref{tilde}). There seems no similar technique and result in existing literature.

Let $\tilde{V}_h^p({\mathcal T}_h)\subset H^1(\Omega)$ denote the continuous piecewise $p$-order finite element space associated with the partition ${\mathcal T}_h$.
For a given function $v\in V_h^p({\mathcal T}_h)$ ($\nsubseteq H^1(\Omega)$), we want to construct a correction function $\tilde{v}\in \tilde{V}_h^p({\mathcal T}_h)$,
which should satisfy the estimates stated in Lemma \ref{tilde}.

Let $v\in V_h^p({\mathcal T}_h)$. For each element $\Omega_k$, we set $v|_{\Omega_k} = v_k$, which denotes the restriction of $v$ on the element $\Omega_k$. We need only to define a suitable
correction function $\tilde{v}_k$ of $v_k$ for each $\Omega_k$. After it is done, we then define the desired function $\tilde{v}$ such that $\tilde{v}|_{\Omega_k}=\tilde{v}_k$.
For ease of understanding, we want to describe the basic idea for defining such function $\tilde{v}_k$. Consider the standard decomposition
$$ v_k=v_k^0+v_k^\partial, $$
where $v_k^\partial|_{\partial\Omega_k}=v_k|_{\partial\Omega_k}$ and $v_k^{\partial}\in V^p_h(\Omega_k)$ is the discrete harmonic extension of $v_k|_{\partial\Omega_k}$ into $\Omega_k$.
It is easy to see that $v_k^0|_{\partial\Omega_k}=v_k|_{\partial\Omega_k}-v_k^\partial|_{\partial\Omega_k}=0$, which can be naturally extended into $\Omega$. However, in general we have $v^{\partial}_k|_{\gamma_{kj}}\not=v^{\partial}_j|_{\gamma_{kj}}$, where $\gamma_{kj}=\partial\Omega_k\cap\partial\Omega_j$ is an element edge.

Since we require that the desired function $\tilde{v}\in H^1(\Omega)$, we need to define a correction $\tilde{v}_k^\partial$ of $v^{\partial}_k$ in a special manner such that $\tilde{v}^{\partial}_k|_{\gamma_{kj}}=\tilde{v}^{\partial}_j|_{\gamma_{kj}}$. After it is done, we naturally define
$$ \tilde{v}_k=v_k^0+\tilde{v}_k^\partial, $$
where $\tilde{v}_k^{\partial}\in V^p_h(\Omega_k)$ is the discrete harmonic extension of $\tilde{v}_k^\partial|_{\partial\Omega_k}$ into $\Omega_k$.

In the following we give a definition of $\tilde{v}_k^\partial|_{\partial\Omega_k}$. Let $e$ denote an edge of $\partial\Omega_k$. When $e = \partial\Omega_k\cap\partial\Omega$, we simply define $\tilde{v}_k^\partial|_e=v^\partial_k|_e$. If $e = \partial\Omega_k\cap\partial\Omega_j$, we define $\tilde{v}_k^\partial|_e$ as follows.

As in the beginning of Subsection \ref{51}, we can define the spaces ${\mathcal P}_1^*$ and ${\mathcal P}'_{p}$ on the edge $e$ by the standard scaling technique.
Then we have the decomposition
\begin{equation*}
v_k^\partial|_e
=v_{k1}^\partial+v_{k0}^\partial,
\end{equation*}
where $v_{k1}^\partial\in{\mathcal P}_1^*$ and $v_{k0}^\partial\in{\mathcal P}'_{p}$. Let $\{\phi_1^{*e},\phi_2^{*e}\}$ denote the two basis functions of ${\mathcal P}_1^*$,
and let $\text{v}_1$ and $\text{v}_2$ denote the two endpoints of the edge $e$. It is easy to see that $v_{k1}^\partial$ can be written as
\begin{equation*}
v_{k1}^\partial
= v_k(\text{v}_1)\phi_1^{*e}+v_k(\text{v}_2)\phi_2^{*e}.
\end{equation*}
%where $\text{v}_1$ and $\text{v}_2$ denotes the two endpoints of the edge $e$.

Set
$$\Lambda_{\text{v}_i}=\{r,~ \Omega_r~\text{contains}~\text{v}_i~ \text{as one of its vertices}\}\quad (i=1,2), $$
and let $n_{\text{v}_i}$ denote the number of all the elements that contain $\text{v}_i$ as their common vertex, namely, the dimension of set $\Lambda_{\text{v}_i}$.
For $e=\partial\Omega_k\cap\partial\Omega_j$, define
\begin{equation} \label{tilde1}
\tilde{v}_{k1}^\partial|_e
= \frac{1}{n_{\text{v}_1}}\sum_{r\in\Lambda_{\text{v}_1}}\!\! v_r(\text{v}_1) \,\phi_1^{*e} +\frac{1}{n_{\text{v}_2}}\sum_{r\in\Lambda_{\text{v}_2}}\!\! v_r(\text{v}_2) \, \phi_2^{*e}
\end{equation}
and
\begin{equation} \label{tilde2}
\tilde{v}_{k0}^\partial|_e=\frac{1}{2}(v_{k0}^\partial+v_{j0}^\partial).
\end{equation}

Now we define $\tilde{v}_k^\partial|_e=\tilde{v}_{k1}^\partial|_e+\tilde{v}_{k0}^\partial|_e$ for each $e\subset \partial\Omega_k$, and let
$\tilde{v}_k^\partial\in V_h^p(\Omega_k)$ be the discrete harmonic extension of $\tilde{v}_k^\partial|_{\partial\Omega_k}$.
From the definition of $\tilde{v}^{\partial}_{k1}$, we know that $\tilde{v}^{\partial}_k|_{\gamma_{kj}}=\tilde{v}^{\partial}_j|_{\gamma_{kj}}$.
Thus we can define $\tilde{v}_k=\tilde{v}_k^\partial+v_k^0$. It is clear that $\tilde{v}_k|_{\gamma_{kj}}=\tilde{v}_j|_{\gamma_{kj}}$.

Finally we define $\tilde{v}$ by $\tilde{v}|_{\Omega_k}=\tilde{v}_k$ and we have $\tilde{v}\in \tilde{V}_h^p({\mathcal T}_h)$.

\begin{lemma} \label{tilde}
For $v\in V_h^p({\mathcal T}_h)$, let $\tilde{v}\in \tilde{V}_h^p({\mathcal T}_h)$ be defined above. Then we have
\begin{equation} \label{tildeeq}
(\sum_{k=1}^N|v-\tilde{v}|^2_{1,\Omega_k}+h^{-2}||v-\tilde{v}||^2_{0,\Omega})^{{1\over 2}}
\lesssim h^{-\frac{1}{2}} p^\frac{1}{2} \big(\sum_{\gamma_{kj}} ||[v]||^2_{0,\gamma_{kj}}\big)^\frac{1}{2}.
\end{equation}
\end{lemma}
\paragraph*{Proof.} Notice that $v_k=v_k^0+v_k^\partial$ and $\tilde{v}_k=v_k^0+\tilde{v}_k^\partial$, and
using the stability of the discrete harmonic extension, we deduce that
\begin{equation} \label{til1}
\begin{split}
\sum_{k=1}^N|v-\tilde{v}|^2_{1,\Omega_k}
& = \sum_{k=1}^N |v_k-\tilde{v}_k|^2_{1,\Omega_k}
= \sum_{k=1}^N |v_k^\partial-\tilde{v}_k^\partial|^2_{1,\Omega_k} \\
& \lesssim \sum_{k=1}^N |v_k^\partial-\tilde{v}_k^\partial|^2_{\frac{1}{2},\partial\Omega_k}
\lesssim h^{-1}p\sum_{k=1}^N ||v_k^\partial-\tilde{v}_k^\partial||^2_{0,\partial\Omega_k}
\end{split}
\end{equation}
and
\begin{flalign} \label{til2}
||v-\tilde{v}||^2_{0,\Omega}
& = \sum_{k=1}^N ||v_k-\tilde{v}_k||^2_{0,\Omega_k}
= \sum_{k=1}^N ||v_k^\partial-\tilde{v}_k^\partial||^2_{0,\Omega_k} \\
& \lesssim \sum_{k=1}^N ( h^2|v_k^\partial-\tilde{v}_k^\partial|^2_{1,\Omega_k} +h||v_k^\partial-\tilde{v}_k^\partial||^2_{0,\partial\Omega_k} )
\lesssim %\sum_{k=1}^N (hp+h) ||v_k^\partial-\tilde{v}_k^\partial||^2_{0,\partial\Omega_k}
hp \sum_{k=1}^N ||v_k^\partial-\tilde{v}_k^\partial||^2_{0,\partial\Omega_k}. \nonumber
\end{flalign}

It suffices to give an estimate of $ \sum\limits_{k=1}^N ||v_k^\partial-\tilde{v}_k^\partial||^2_{0,\partial\Omega_k}$.

Let $e$ be an edge on $\partial\Omega_k$. When $e = \partial\Omega_k\cap\partial\Omega$, we have $ (v_k^\partial-\tilde{v}_k^\partial)|_e
=0$, which implies that $||v_k^\partial-\tilde{v}_k^\partial||_{0,e}=0$.

If $e = \partial\Omega_k\cap\partial\Omega_j=\gamma_{kj}$, we have
$(v_k^\partial-\tilde{v}_k^\partial)|_{\gamma_{kj}}
= (v_{k1}^\partial-\tilde{v}_{k1}^\partial)+(v_{k0}^\partial-\tilde{v}_{k0}^\partial)$. So we get
\begin{equation} \label{zong}
||v_k^\partial-\tilde{v}_k^\partial||_{0,\gamma_{kj}}
\leq ||v_{k1}^\partial-\tilde{v}_{k1}^\partial||_{0,\gamma_{kj}} +||v_{k0}^\partial-\tilde{v}_{k0}^\partial||_{0,\gamma_{kj}}.
\end{equation}

Let v$_1$ and v$_2$ denote two endpoints of $\gamma_{kj}$, and let $\Lambda_{\text{v}_i}$ $(i=1,2)$ be the sets defined before (\ref{tilde1}).
It follows, from (\ref{tilde1}) and (\ref{tilde2}), that
\begin{equation} \label{fir}
v_{k1}^\partial-\tilde{v}_{k1}^\partial
= \frac{1}{n_{\text{v}_1}}\sum_{r\in\Lambda_{\text{v}_1}} (v_k-v_r)(\text{v}_1) \,\phi_1^{*e} +\frac{1}{n_{\text{v}_2}}\sum_{r\in\Lambda_{\text{v}_2}} (v_k-v_r)(\text{v}_2) \, \phi_2^{*e}
\end{equation}
and
\begin{equation*}
v_{k0}^\partial-\tilde{v}_{k0}^\partial
= \frac{1}{2}(v_{k0}^\partial-v_{j0}^\partial).
\end{equation*}

Notice that $(v_k^\partial-v_j^\partial)|_{\gamma_{kj}} = (v_{k1}^\partial-v_{j1}^\partial)+(v_{k0}^\partial-v_{j0}^\partial)$, where $(v_{k1}^\partial-v_{j1}^\partial)\in{\mathcal P}_1^*$ and $(v_{k0}^\partial-v_{j0}^\partial)\in{\mathcal P}'_{p}$. It is easy to see that
\begin{equation}
||v_{k0}^\partial-\tilde{v}_{k0}^\partial||_{0,\gamma_{kj}}
= \frac{1}{2}||v_{k0}^\partial-v_{j0}^\partial||_{0,\gamma_{kj}}
\leq \frac{1}{2}||v_k^\partial-v_j^\partial||_{0,\gamma_{kj}}
= \frac{1}{2}||v_k-v_j||_{0,\gamma_{kj}}.\label{4.inequality-new1}
\end{equation}

Clearly, we have $k\in\Lambda_{\text{v}_i}$ $(i=1,2)$. Define
$$\Lambda_{\text{v}_i}^{k,1}=\{r,~ r\in\Lambda_{\text{v}_i}, r\ne k, \text{ and } \gamma_{kr}=\partial\Omega_k\cap\partial\Omega_r \text{ is an edge} \}\quad (i=1,2) $$
and
$$\Lambda_{\text{v}_i}^{k,2}=\{r,~ r\in\Lambda_{\text{v}_i}, r\ne k, \text{ and } \partial\Omega_r\cap\partial\Omega_j \text{ is the vertex v}_i  \}\quad (i=1,2). $$
It is clear that $\Lambda_{\text{v}_i}=\{k\}\cup\Lambda_{\text{v}_i}^{k,1}\cup\Lambda_{\text{v}_i}^{k,2}$ and $\Lambda_{\text{v}_i}^{k,1}\cap\Lambda_{\text{v}_i}^{k,2}=\emptyset$.
Then, for the first item on the right side of (\ref{fir}), we have
\begin{eqnarray}
\frac{1}{n_{\text{v}_1}}\sum_{r\in\Lambda_{\text{v}_1}} (v_k-v_r)(\text{v}_1) \,\phi_1^{*e}
&=& \frac{1}{n_{\text{v}_1}}\sum_{r\in\Lambda_{\text{v}_1}^{k,1}} (v_k-v_r)(\text{v}_1) \,\phi_1^{*e}\cr
&+& \frac{1}{n_{\text{v}_1}}\sum_{r\in\Lambda_{\text{v}_1}^{k,2}} (v_k-v_r)(\text{v}_1) \,\phi_1^{*e}. \label{4.estimate-new0}
\end{eqnarray}

We first give an estimate of $||(v_k-v_r)(\text{v}_1) \,\phi_1^{*e}||_{0,e}$ for $r\in\Lambda_{\text{v}_1}^{k,1}$.
Let $\{\phi_1^{*\gamma_{kr}}, \phi_2^{*\gamma_{kr}}\}$ denote the two basis functions
of ${\mathcal P}_1^*$ associated with the edge $\gamma_{kr}$. We already assume that the partition ${\mathcal T}_h$ is quasi-uniform, which yields
\begin{equation*}
||\phi_1^{*e}||_{0,e}
\lesssim ||\phi_1^{*\gamma_{kr}}||_{0,\gamma_{kr}}
\text{ and }
||\phi_2^{*e}||_{0,e}
\lesssim ||\phi_2^{*\gamma_{kr}}||_{0,\gamma_{kr}}.
\end{equation*}
This, together with Lemma \ref{star} (replacing $q$ by $p$), leads to
\begin{equation}
\begin{split}
||(v_k-v_r)(\text{v}_1) \,\phi_1^{*e}||_{0,e}
& \lesssim |v_k-v_r|(\text{v}_1) \,||\phi_1^{*\gamma_{kr}}||_{0,\gamma_{kr}}\\
&\leq \sqrt{\frac{p+1}{p}} ||v^\partial_{k1}-v^\partial_{r1}||_{0,\gamma_{kr}}
\lesssim ||v^\partial_k-v^\partial_r||_{0,\gamma_{kr}}\\
&= ||v_k-v_r||_{0,\gamma_{kr}}
= ||[v]||_{0,\gamma_{kr}}\quad\quad(r\in \Lambda_{\text{v}_1}^{k,1}). \label{4.estimate-new1}
\end{split}
\end{equation}

 Next we estimate $||(v_k-v_r)(\text{v}_1) \,\phi_1^{*e}||_{0,e}$ for $r\in\Lambda_{\text{v}_1}^{k,2}$.
Without loss of generality, we assume that there exists some element $\Omega_l$ such that $\partial\Omega_k\cap\partial\Omega_l=\gamma_{kl}$
and $\partial\Omega_l\cap\partial\Omega_r=\gamma_{lr}$ are two (different) edges that have the common vertex $\text{v}_1$. Then, by the triangle inequality, we have
$$ ||(v_k-v_r)(\text{v}_1) \,\phi_1^{*e}||_{0,e}\leq ||(v_k-v_l)(\text{v}_1) \,\phi_1^{*e}||_{0,e}+||(v_l-v_r)(\text{v}_1) \,\phi_1^{*e}||_{0,e}.$$
We can estimate the two terms on the right side of the above inequality like (\ref{4.estimate-new1}), and we obtain
\begin{equation} ||(v_k-v_r)(\text{v}_1) \,\phi_1^{*e}||_{0,e}\lesssim ||[v]||_{0,\gamma_{kl}}+||[v]||_{0,\gamma_{lr}}\quad\quad(r\in \Lambda_{\text{v}_1}^{k,2}).\label{4.estimate-new2}
\end{equation}
Let ${\mathcal E}_{\text{v}_1}$ denote the set of all the edges that have $\text{v}_1$ as their common endpoint. Substituting (\ref{4.estimate-new1}) and (\ref{4.estimate-new2}) into
(\ref{4.estimate-new0}), yields
\begin{equation}
||\frac{1}{n_{\text{v}_1}}\!\sum_{r\in\Lambda_{\text{v}_1}}\! (v_k-v_r)(\text{v}_1) \,\phi_1^{*e}||_{0,e}
\!\lesssim \sum\limits_{e\in{\mathcal E}_{\text{v}_1}} ||[v]||_{0,e}. \label{4.estimate-new3}
\end{equation}

In an analogous way with (\ref{4.estimate-new3}), we can verify that
\begin{equation*}
||\frac{1}{n_{\text{v}_2}}\!\sum_{r\in\Lambda_{\text{v}_2}}\! (v_k-v_r)(\text{v}_2) \,\phi_2^{*e}||_{0,e}
\!\lesssim \sum\limits_{e\in{\mathcal E}_{\text{v}_2}} ||[v]||_{0,e}.
\end{equation*}
Here ${\mathcal E}_{\text{v}_2}$ denotes the set of all the edges that have $\text{v}_2$ as their common endpoint.
Plugging this and (\ref{4.estimate-new3}) in (\ref{fir}), leads to
$$ ||v_{k1}^\partial-\tilde{v}_{k1}^\partial||_{0,\gamma_{kj}}\!\lesssim \sum\limits_{e\in{\mathcal E}_{\text{v}_1}} ||[v]||_{0,e}+\sum\limits_{e\in{\mathcal E}_{\text{v}_2}} ||[v]||_{0,e}, $$
which, together with (\ref{4.inequality-new1}) and (\ref{zong}), gives
\begin{equation*}
\begin{split}
||v_k^\partial-\tilde{v}_k^\partial||_{0,\gamma_{kj}}
& \lesssim \sum\limits_{e\in{\mathcal E}_{\text{v}_1}} ||[v]||_{0,e}+\sum\limits_{e\in{\mathcal E}_{\text{v}_2}} ||[v]||_{0,e}.
\end{split}
\end{equation*}
Hence, we get
$$
\sum_{k=1}^N ||v_k^\partial-\tilde{v}_k^\partial||^2_{0,\partial\Omega_k}
= \sum_{\gamma_{kj}} ||v_k^\partial-\tilde{v}_k^\partial||^2_{0,\gamma_{kj}}
\lesssim  \sum_{\gamma_{kj}} ||[v]||^2_{0,\gamma_{kj}}.
$$
Finally, substituting this inequality into (\ref{til1}) and (\ref{til2}), we obtain
\begin{equation*}
\sum_{k=1}^N|v-\tilde{v}|^2_{1,\Omega_k}
\lesssim h^{-1}p \sum_{\gamma_{kj}} ||[v]||^2_{0,\gamma_{kj}}
\text{ and }\
||v-\tilde{v}||^2_{0,\Omega}
\lesssim hp \sum_{\gamma_{kj}} ||[v]||^2_{0,\gamma_{kj}}.
\end{equation*}
The estimate (\ref{tildeeq}) is a direct consequence of the above inequalities.
$\hfill \Box$

In the following we want to build a Poincare-type inequality for the functions in $V_h^p({\mathcal T}_h)$ by Lemma \ref{tilde}.
%By Lemma \ref{tilde} and Lemma \ref{h}, we can derive a $Poincar\acute{e}$-type inequality.\\
\begin{lemma}\label{dual} Let Assumption 1-Assumption 4 be satisfied. Assume that, for some $\lambda_h\in W_h^q(\gamma)$, the function $v\in V_h^p({\mathcal T}_h)$ satisfy
\begin{equation} \label{bb1}
a^{(k)}(v, w)= \langle \pm\lambda_h, \overline{w}\rangle_{\partial\Omega_k \backslash\partial\Omega}\quad (k=1,2,\cdots,N),~~ \forall\,w\in V_h^p({\mathcal T}_h).
\end{equation}
Then
\begin{equation}
||v||^2_{0,\Omega}
\lesssim h^{-1} \sum_{\gamma_{kj}} ||[v]||^2_{0,\gamma_{kj}}.\label{5.Poincare}
\end{equation}
\end{lemma}
\paragraph*{Proof.} A standard technique to estimate $L^2$ norm of a function is the introduction of a suitable dual problem (see, for example, \cite{Hiptmair2011}). Consider the dual problem
%In the dual problem (\ref{adjo1}), we choose $\tilde{f} = v$. Then (\ref{adjo1}) becomes %{\mathbf \kappa}^2
\begin{equation} \label{adjo}
\left\{
\begin{aligned}
& -\Delta\phi - \kappa^2 \phi = v
& \text{in}\ \Omega,\\
& \frac{\partial \phi}{\partial n} - i\kappa \phi = 0
& \text{on}\ \partial\Omega.
\end{aligned}
\right.
\end{equation}
Let $\phi\in H^1(\Omega)$ and $\phi_h\in \tilde{V}_h^p({\mathcal T}_h)$ denote its weak solution and $p$-order finite element solution, which are defined respectively by
\begin{equation} \label{fem1}
(\nabla \phi,\overline{\nabla \psi})_\Omega-(\kappa^2\phi,\overline{\psi})_\Omega-i\langle \kappa\phi,\overline{\psi} \rangle_{\partial\Omega}
= (v,\overline{\psi})_\Omega, \quad \forall \psi\in H^1(\Omega)
\end{equation}
and
\begin{equation} \label{fem2}
(\nabla \phi_h,\overline{\nabla \psi}_h)_\Omega-(\kappa^2\phi_h,\overline{\psi}_h)_\Omega-i\langle\kappa\phi_h,\overline{\psi}_h \rangle_{\partial\Omega}
= (v,\overline{\psi}_h)_\Omega, \quad \forall \psi_h\in \tilde{V}_h^p({\mathcal T}_h).
\end{equation}

Using (\ref{adjo}) and Green's formula, we obtain
\begin{flalign}
||v||^2_{0,\Omega}
& = \sum_{k=1}^N (v,\overline{v})_{\Omega_k}
= \sum_{k=1}^N \big( (v,\overline{-\Delta\phi})_{\Omega_k}-(v,\overline{\kappa^2\phi})_{\Omega_k} \big) \nonumber\\
& = \sum_{k=1}^N \big( (\nabla v,\overline{\nabla\phi})_{\Omega_k} -\langle v,\overline{\nabla\phi\!\cdot\!n}\rangle_{\partial\Omega_k} -(\kappa^2 v,\overline{\phi})_{\Omega_k} \big) \nonumber\\
& = \sum_{k=1}^N (\nabla v,\overline{\nabla\phi})_{\Omega_k}-\sum_{k=1}^N (\kappa^2 v,\overline{\phi})_{\Omega_k}-\sum_{\gamma_{kj}} \langle[v],\overline{\nabla\phi\!\cdot\!n}\rangle_{\gamma_{kj}}-\langle v,\overline{i\kappa\phi}\rangle_{\partial\Omega} \nonumber\\
& = \sum_{k=1}^N (\nabla v,\overline{\nabla\phi_h})_{\Omega_k}-\sum_{k=1}^N(\kappa^2 v,\overline{\phi_h})_{\Omega_k} -\sum_{\gamma_{kj}} \langle[v],\overline{\nabla\phi\!\cdot\!n}\rangle_{\gamma_{kj}}
+i\langle\kappa v,\overline{\phi}\rangle_{\partial\Omega} \nonumber\\
& +\sum_{k=1}^N (\nabla v,\overline{\nabla(\phi-\phi_h)})_{\Omega_k}-\sum_{k=1}^N(\kappa^2v,\overline{\phi-\phi_h})_{\Omega_k}.\label{4.equality-new1}
\end{flalign}
In the last equality, we introduced the finite element function $\phi_h$ since (\ref{bb1}) holds only for finite element function $w$ (if $v\in H^1(\Omega_k)$ satisfies the equation (\ref{bb1}) for any $w\in H^1(\Omega_k)$,
then the proof is trivial).

Letting $w = \phi_h$ in (\ref{bb1}) and summing the resulting equality over $k$, and using the fact that $\phi_h$ is continuous across the inner edges, gives
\begin{equation*}
\begin{split}
\sum_{k=1}^N a^{(k)}(v, \phi_h)
= \sum_{k=1}^N \langle \pm\lambda_h, \overline{\phi_h}\rangle_{\partial\Omega_k \backslash\partial\Omega}
=0,
\end{split}
\end{equation*}
which implies that
\begin{equation*}
\sum_{k=1}^N (\nabla v, \overline{\nabla \phi}_h)_{\Omega_k}-\sum_{k=1}^N(\kappa^2v, \overline{\phi}_h)_{\Omega_k}
= -i\rho\sum_{\gamma_{kj}}\langle [v], \overline{\phi}_h\rangle_{\gamma_{kj}}-i\langle\kappa v, \overline{\phi}_h\rangle_{\partial\Omega}.
\end{equation*}
This, together with (\ref{4.equality-new1}), leads to
\begin{flalign} \label{s40}
||v||^2_{0,\Omega}
& = -i\rho\sum_{\gamma_{kj}}\langle [v], \overline{\phi_h}\rangle_{\gamma_{kj}}-i\langle\kappa v, \overline{\phi}_h\rangle_{\partial\Omega} -\sum_{\gamma_{kj}} \langle[v],\overline{\nabla\phi\!\cdot\!n}\rangle_{\gamma_{kj}}+i\langle\kappa v, \overline{\phi}\rangle_{\partial\Omega} \nonumber\\
& +\sum_{k=1}^N (\nabla v,\overline{\nabla(\phi-\phi_h)})_{\Omega_k}-\sum_{k=1}^N (\kappa^2v,\overline{\phi-\phi_h})_{\Omega_k} \nonumber\\
& = \sum_{k=1}^N (\nabla v,\overline{\nabla(\phi-\phi_h)})_{\Omega_k}-\sum_{k=1}^N(\kappa^2v,\overline{\phi-\phi_h})_{\Omega_k}+i\langle\kappa v, \overline{\phi-\phi_h}\rangle_{\partial\Omega} \nonumber\\
& -\sum_{\gamma_{kj}} \langle[v],\overline{\nabla\phi\!\cdot\!n}\rangle_{\gamma_{kj}} -i\rho\sum_{\gamma_{kj}}\langle[v],\overline{\phi}_h\rangle_{\gamma_{kj}}.
\end{flalign}
If we directly estimate the terms containing the error $\phi-\phi_h$, we cannot build the inequality (\ref{5.Poincare}) unless a stronger assumption on the mesh size $h$ is made. Because of this,
we have to introduce a globally continuous finite element ``approximation" of $v$ such that the energy orthogonality of $\phi-\phi_h$ can be used.

For $v\in V_h^p({\mathcal T}_h)$, we construct $\tilde{v}\in \tilde{V}_h^p({\mathcal T}_h)$ as in Lemma \ref{tilde}. For ease of notation, set
\begin{equation*}
R
=\sum_{k=1}^N\! (\nabla (v-\tilde{v}),\overline{\nabla(\phi-\phi_h)})_{\Omega_k}\!-\!\sum_{k=1}^N\! (\kappa^2(v-\tilde{v}),\overline{\phi-\phi_h})_{\Omega_k}\!+\!i\langle \kappa(v-\tilde{v}),\overline{\phi-\phi_h}\rangle_{\partial\Omega}.
\end{equation*}
Then (\ref{s40}) can be written as
\begin{flalign}  \label{s41}
||v||^2_{0,\Omega}
& = R\!+\!\sum_{k=1}^N (\nabla \tilde{v},\overline{\nabla(\phi-\phi_h)})_{\Omega_k}\!-\!\sum_{k=1}^N (\kappa^2\tilde{v},\overline{\phi-\phi_h})_{\Omega_k}\!+\!i\langle\kappa \tilde{v},\overline{\phi-\phi_h}\rangle_{\partial\Omega} \nonumber\\
&  -\sum_{\gamma_{kj}} \langle[v],\!\overline{\nabla\phi\!\cdot\!n}\rangle_{\gamma_{kj}} -i\rho\sum_{\gamma_{kj}}\langle[v],\overline{\phi}_h\rangle_{\gamma_{kj}}.
\end{flalign}

Choosing $\psi = \tilde{v}$ in (\ref{fem1}) and $\psi_h = \tilde{v}$ in (\ref{fem2}), we get the difference
\begin{equation} \label{fem31}
(\nabla (\phi-\phi_h),\overline{\nabla \tilde{v}})_\Omega-(\kappa^2(\phi-\phi_h),\overline{\tilde{v}})_\Omega-i\langle\kappa(\phi-\phi_h),\overline{\tilde{v}} \rangle_{\partial\Omega}
= 0,
\end{equation}
which is called as the energy orthogonality of $\phi-\phi_h$. The complex conjugation of (\ref{fem31}) becomes
\begin{equation} \label{fem41}
(\nabla \tilde{v},\overline{\nabla (\phi-\phi_h)})_\Omega-(\kappa^2\tilde{v},\overline{\phi-\phi_h})_\Omega+i\langle\kappa \tilde{v},\overline{\phi-\phi_h} \rangle_{\partial\Omega}
= 0.
\end{equation}
Substituting (\ref{fem41}) into (\ref{s41}), we obtain
\begin{equation} \label{s42}
||v||^2_{0,\Omega}
= R-\sum_{\gamma_{kj}} \langle[v],\overline{\nabla\phi\!\cdot\!n}\rangle_{\gamma_{kj}} -i\rho\sum_{\gamma_{kj}}\langle[v],\overline{\phi}_h\rangle_{\gamma_{kj}}.
\end{equation}

Let $M=\sum\limits_{\gamma_{kj}} ||[v]||^2_{0,\gamma_{kj}}$. Using Cauchy-Schwarz inequality to the sums on the right side of (\ref{s42}), yields
\begin{flalign}  \label{s51}
||v||^2_{0,\Omega}
& \leq |R|+\sum_{\gamma_{kj}} ||[v]||_{0,\gamma_{kj}} ||\nabla\phi\!\cdot\!n||_{0,\gamma_{kj}} +\rho\sum_{\gamma_{kj}}||[v]||_{0,\gamma_{kj}} ||\phi_h||_{0,\gamma_{kj}} \nonumber\\
& \leq |R|+(\sum_{\gamma_{kj}} ||\nabla\phi\!\cdot\!n||^2_{0,\gamma_{kj}})^\frac{1}{2}M^\frac{1}{2} +\rho (\sum_{\gamma_{kj}} ||\phi_h||^2_{0,\gamma_{kj}})^\frac{1}{2}M^\frac{1}{2}.
\end{flalign}
It suffices to estimate $|R|$. It is easy to see that
\begin{flalign} \label{r}
|R|
&\lesssim\sum_{k=1}^N ||\nabla (v-\tilde{v})||_{0,\Omega_k} ||\nabla(\phi-\phi_h)||_{0,\Omega_k}+\sum_{k=1}^N \omega^2||v-\tilde{v}||_{0,\Omega_k} ||\phi-\phi_h||_{0,\Omega_k} \nonumber\\
& +\omega||v-\tilde{v}||_{0,\partial\Omega} ||\phi-\phi_h||_{0,\partial\Omega} \nonumber\\
&\lesssim(\sum_{k=1}^N||\nabla (v-\tilde{v})||^2_{0,\Omega_k})^{{1\over 2}}||\nabla(\phi-\phi_h)||_{0,\Omega} +\omega^2||v-\tilde{v}||_{0,\Omega} ||\phi-\phi_h||_{0,\Omega} \nonumber\\
& +\omega||v-\tilde{v}||_{0,\partial\Omega} ||\phi-\phi_h||_{0,\partial\Omega}.
\end{flalign}
It follows by Assumption 4 that
\begin{equation*}
||\nabla(\phi-\phi_h)||_{0,\Omega} + \omega||\phi-\phi_h||_{0,\Omega}
\lesssim p^{-{1\over 2}} ||v||_{0,\Omega}.
\end{equation*}
Then, by the $\varepsilon$-inequality ($\varepsilon=\omega^{-\frac{1}{2}}<1$), we have
\begin{equation*}
||\phi-\phi_h||_{0,\partial\Omega}
\lesssim \omega^{-\frac{1}{2}} ||\nabla(\phi-\phi_h)||_{0,\Omega} \!+\! \omega^\frac{1}{2} ||\phi-\phi_h||_{0,\Omega}
\lesssim \omega^{-\frac{1}{2}}p^{-{1\over 2}} ||v||_{0,\Omega}.
\end{equation*}
Moreover, from Lemma \ref{tilde}, we have
\begin{equation*}
(\sum_{k=1}^N||\nabla(v-\tilde{v})||^2_{0,\Omega_k} + h^{-2}||v-\tilde{v}||^2_{0,\Omega})^{{1\over 2}}
\lesssim h^{-\frac{1}{2}}p^\frac{1}{2} M^\frac{1}{2}.
\end{equation*}
Furthermore, by the trace inequality, we get
\begin{flalign*}
||v-\tilde{v}||^2_{0,\partial\Omega}
& = \sum_{k=1}^N ||v-\tilde{v}||^2_{0,\partial\Omega_k \cap\partial\Omega}
%\leq \sum_{k=1}^N ||v-\tilde{v}||^2_{0,\partial\Omega_k} \\
\lesssim \sum_{k=1}^N \big( h||\nabla(v-\tilde{v})||^2_{0,\Omega_k} +h^{-1}||v-\tilde{v}||^2_{0,\Omega_k} \big) \nonumber\\
& = h||\nabla(v-\tilde{v})||^2_{0,\Omega}+h^{-1}||v-\tilde{v}||^2_{0,\Omega}
\lesssim p M.
\end{flalign*}
Hence, substituting the above estimates into (\ref{r}), we obtain
\begin{equation} \label{r1}
|R|
\lesssim (1+\omega h+
+\omega^\frac{1}{2}h^{\frac{1}{2}})h^{-\frac{1}{2}}M^\frac{1}{2} ||v||_{0,\Omega}.
\end{equation}

On the other hand, by Assumption 1 and Theorem 1 of \cite{Brown2016} the function $\phi$ satisfies
$$ |\phi|_{1,\Omega}+\omega\|\phi\|_{0,\Omega}\lesssim\|v\|_{0,\Omega}. $$
As in the proof of Lemma 3.3 of \cite{DuWu2015}, we can further verify that $|\phi|_{2,\Omega}\lesssim\omega\|v\|_{0,\Omega}$.
Thus we have the stability
\begin{equation*}
\omega^{-1}|\phi|_{2,\Omega}+|\phi|_{1,\Omega}+\omega\|\phi\|_{0,\Omega}
\lesssim ||v||_{0,\Omega}.
\end{equation*}
Then, by the trace inequality, we get
\begin{equation} \label{s71}
\begin{split}
\sum_{\gamma_{kj}} ||\nabla\phi\!\cdot\!n||^2_{0,\gamma_{kj}}
& \leq \sum_{k=1}^N ||\nabla\phi\!\cdot\!n||^2_{0,\partial\Omega_k}
\lesssim \sum_{k=1}^N (h|\phi|^2_{2,\Omega_k} +h^{-1}|\phi|^2_{1,\Omega_k}) \cr
& = h|\phi|^2_{2,\Omega} +h^{-1}|\phi|^2_{1,\Omega}
\lesssim (h^2\omega^2+1)h^{-1} ||v||^2_{0,\Omega}.
\end{split}
\end{equation}
and (using Assumption 4 again)
\begin{equation} \label{s72}
\begin{split}
\sum_{\gamma_{kj}} ||\phi_h||^2_{0,\gamma_{kj}}
&\leq \sum_{k=1}^N ||\phi_h||^2_{0,\partial\Omega_k}
\lesssim \sum_{k=1}^N (h|\phi_h|^2_{1,\Omega_k} +h^{-1}\|\phi_h\|^2_{0,\Omega_k}) \\
&\lesssim h|\phi|^2_{1,\Omega} +h^{-1}\|\phi\|^2_{0,\Omega}+(h+h^{-1}\omega^{-2})p^{-1}\|v\|^2_{0,\Omega}\\
&\lesssim(h^2+\omega^{-2})h^{-1} ||v||^2_{0,\Omega}.
\end{split}
\end{equation}

Substituting the inequalities (\ref{r1}), (\ref{s71}) and (\ref{s72}) into (\ref{s51}) and using Assumption 2 and Assumption 3, yields
\begin{equation*}
||v||^2_{0,\Omega}
\lesssim h^{-\frac{1}{2}} \cdot M^\frac{1}{2} || v||_{0,\Omega}.
\end{equation*}
Finally, we obtain the desired inequality (\ref{5.Poincare}).
$\hfill \Box$
\begin{remark} The inequality (\ref{5.Poincare}) can be viewed as an extension of the Poincare inequality held for plane wave functions to the piecewise polynomial functions in $V_h^p({\mathcal T}_h)$.
Comparing the inequality (\ref{5.Poincare}) with the Poincare-type inequality given by Lemma 3.7 of \cite{Hiptmair2011} for the plane wave functions, we find that the right sides of the two inequalities contain
the same term $h^{-1}\sum_{\gamma_{kj}} ||[v]||^2_{0,\gamma_{kj}}$, and (\ref{5.Poincare}) is more succinct thanks to the condition (\ref{bb1}) (there are extra terms in the inequality in Lemma 3.7 of \cite{Hiptmair2011}). However, the proof of
(\ref{5.Poincare}), which depends on the estimates (\ref{tildeeq}) and (\ref{5.2new06}), is much more technical than that of the inequality in Lemma 3.7 of \cite{Hiptmair2011}
since the considered functions do not satisfy the homogeneous Helmholtz equation satisfied by the plane wave functions.
\end{remark}

\begin{remark}
As pointed out in Remark 2.1, the proposed method is practical for the case with variable wave numbers, but we have to use Assumption 4 to give the theoretical
analysis of (\ref{5.Poincare}). The main reason is that we do not know whether the estimate (\ref{5.2new06}) proved in \cite{fem} for constant wave numbers is still valid for the case with
variable wave numbers (the condition (\ref{4.assumption4}) that we used is much weaker than (\ref{5.2new06})). We failed to build a similar inequality with (\ref{5.Poincare}) without Assumption 4.
\end{remark}

 In the rest of this paper, we always use $a^{(k)}(\cdot,\cdot)$ to denote the local sesquilinear form defined in Subsection 2.2. For $v\in \prod_{k=1}^N H^1(\Omega_k)$, define
% $$ \||v\||_k=(\|\nabla v\|^2_{0,\Omega_k}+\omega^2\|v\|^2_{0,\Omega_k})^{{1\over 2}},\quad k=1,\cdots,N $$
%and
$$ \||v\||=(\sum\limits_{k=1}^N\|\nabla v\|^2_{0,\Omega_k}+\omega^2\|v\|^2_{0,\Omega})^{{1\over 2}}. $$

\begin{lemma} \label{lamda} Let Assumption 3 be satisfied. Suppose $q\geq 1$ and $p\geq q+2$. Assume that $v\in V_h^p({\mathcal T}_h)$ and $\lambda_h\in W_h^q(\gamma)$ satisfy the relation
\begin{equation} \label{bo}
a^{(k)}(v, w)= \langle \pm\lambda_h, \overline{w}\rangle_{\partial\Omega_k \backslash\partial\Omega}\quad\quad(k=1,2,\cdots,N),~~\forall\,w\in V_h^p({\mathcal T}_h).
\end{equation}
Then the following estimate holds
\begin{equation}
\sum_{\gamma_{kj}} ||\lambda_h||^2_{0,\gamma_{kj}}\lesssim h^{-1}pq\||v\||^2.\label{5.2new00}
\end{equation}
%where $\delta=1$ when $p\geq q+2$ or $\delta=0$ if $p=2q$ ($q\geq 2$).
\end{lemma}

\paragraph*{Proof.}
By the local $inf-sup$ condition given in Theorem \ref{infsup}, there exists a non-zero function $\psi\in V_h^p(\Omega_k)$, which is discrete harmonic in $\Omega_k$,
such that
\begin{equation*}
||\lambda_h||^2_{0,\partial\Omega_k\backslash\partial\Omega}=\|\pm\lambda_h||^2_{0,\partial\Omega_k\backslash\partial\Omega}
\lesssim q \frac{\langle \pm\lambda_h, \overline{\psi} \rangle^2_{\partial\Omega_k\backslash\partial\Omega}} {||\psi||^2_{0,\partial\Omega_k}}.
\end{equation*}
%where $\delta=1$ when $p\geq q+2$ or $\delta=0$ if $p=2q$ ($q\geq 2$).
Then, using (\ref{bo}), Cauchy inequality and Assumption 3, yields
\begin{equation}
\begin{split}
& ||\lambda_h||^2_{0,\partial\Omega_k\backslash\partial\Omega}
\lesssim q \cdot \frac{\big(a^{(k)}(v, \psi)\big)^2}{||\psi||^2_{0,\partial\Omega_k}} \\
& \lesssim q \frac{||\nabla v||^2_{0,\Omega_k} ||\nabla \psi||^2_{0,\Omega_k}+\omega^4||v||^2_{0,\Omega_k} ||\psi||^2_{0,\Omega_k}+\omega^2||v||^2_{0,\partial\Omega_k} ||\psi||^2_{0,\partial\Omega_k}}{||\psi||^2_{0,\partial\Omega_k}}. \label{4.inequality-new8}
\end{split}
\end{equation}
By the stability of discrete harmonic functions, the inverse estimate and Poincare inequality, we deduce that
\begin{equation*}
\begin{split}
& ||\nabla \psi||^2_{0,\Omega_k}
\lesssim |\psi|^2_{\frac{1}{2},\partial\Omega_k}
\lesssim h^{-1}p ||\psi||^2_{0,\partial\Omega_k}\\
& ||\psi||^2_{0,\Omega_k}\lesssim h^2||\nabla \psi||^2_{0,\Omega_k}+h||\psi||^2_{0,\partial\Omega_k}
\lesssim hp\|\psi\|^2_{0,\partial\Omega_k}+h||\psi||^2_{0,\partial\Omega_k}
\lesssim hp||\psi||^2_{0,\partial\Omega_k}.
\end{split}
\end{equation*}
Substituting this into (\ref{4.inequality-new8}), together with the trace inequality and Assumption 2, yields
\begin{equation*}
\begin{split}
||\lambda_h||^2_{0,\partial\Omega_k\backslash\partial\Omega}
& \lesssim q \big( ||\nabla v||^2_{0,\Omega_k}h^{-1}p +\omega^4||v||^2_{0,\Omega_k} hp +\omega^2(h||\nabla v||^2_{0,\Omega_k}+h^{-1}||v||^2_{0,\Omega_k})) \\
& = q\big( (h^{-1}p+\omega^2h)||\nabla v||^2_{0,\Omega_k}+(\omega^4hp+\omega^2h^{-1}) ||v||^2_{0,\Omega_k} \big) \\
&=h^{-1}pq\big((1+{\omega^2h^2\over p})||\nabla v||^2_{0,\Omega_k}+(\omega^2h^2+p^{-1})\omega^2||v||^2_{0,\Omega_k} \big)\\
& \lesssim h^{-1}pq(||\nabla v||^2_{0,\Omega_k}+\omega^2||v||^2_{0,\Omega_k}).
\end{split}
\end{equation*}
Summing up the above inequality over $k$, gives (\ref{5.2new00}).
$\hfill \Box$

By Lemma \ref{dual} and Lemma \ref{lamda}, we can prove a crucial auxiliary result given below, which can be viewed as a {\it jump-controlled} stability estimate.
As we will see, this auxiliary result plays a key role in the proof of Theorem 4.2 and Theorem 4.3.\\
{\bf Proposition 5.1} Assume that $q\geq 1$ and $p\geq q+2$. Let Assumption 1-Assumption 4 be satisfied, and let
%$\omega hp^{-1} \leq C_0$ and $p\geq 1+c_0\log\omega$ for positive constants $c_0$ and $C_0$.
$v\in V_h^p({\mathcal T}_h)$ satisfy
\begin{equation} \label{bb}
a^{(k)}(v, w)= \langle \pm\lambda_h, \overline{w}\rangle_{\partial\Omega_k \backslash\partial\Omega}\quad\quad(k=1,2,\cdots,N),~\forall\,w\in V_h^p({\mathcal T}_h).
\end{equation}
Then the following estimate holds
\begin{equation}
\||v\||^2
\lesssim \omega^2h^{-1} \sum_{\gamma_{kj}} ||[v]||^2_{0,\gamma_{kj}}.\label{4.inequality.new10}
\end{equation}
\paragraph*{Proof.}
Choosing $w = v$ in (\ref{bb}) and summing the resulting equality over $k$, gives
\begin{equation*}
\begin{split}
& \sum_{k=1}^N \big( ||\nabla v||^2_{0,\Omega_k}-(\kappa^2v,v)_{\Omega_k} \pm i\rho ||v||^2_{0,\partial\Omega_k \backslash\partial\Omega}+ i\langle\kappa v,v\rangle_{\partial\Omega_k \cap\partial\Omega} \big) \\
& = \sum_{k=1}^N \langle \pm\lambda_h,\overline{v} \rangle_{\partial\Omega_k\backslash\partial\Omega}
= \sum_{\gamma_{kj}} \langle \lambda_h,[\overline{v}] \rangle_{\gamma_{kj}}.
\end{split}
\end{equation*}
Let $M=\sum\limits_{\gamma_{kj}} ||[v]||^2_{0,\gamma_{kj}}$. Considering the module of the above equality and using Cauchy-Schwarz inequality, yields
% and the imaginary part
\begin{equation*}
|\sum_{k=1}^N||\nabla v||^2_{0,\Omega_k}-(\kappa^2v, v)_{\Omega}|
\leq \sum_{\gamma_{kj}} ||\lambda_h||_{0,\gamma_{kj}} ||[v]||_{0,\gamma_{kj}}
\leq \big(\sum_{\gamma_{kj}} ||\lambda_h||^2_{0,\gamma_{kj}}\big)^\frac{1}{2} M^\frac{1}{2},
\end{equation*}
which implies that
$$ \sum_{k=1}^N||\nabla v||^2_{0,\Omega_k}\lesssim \omega^2|| v||^2_{0,\Omega}+\big(\sum_{\gamma_{kj}} ||\lambda_h||^2_{0,\gamma_{kj}}\big)^\frac{1}{2} M^\frac{1}{2}. $$
Combining this inequality with (\ref{5.Poincare}) of Lemma \ref{dual}, leads to
\begin{equation}  \label{b4}
\begin{split}
\sum_{k=1}^N||\nabla v||^2_{0,\Omega_k}+\omega^2|| v||^2_{0,\Omega}
&\lesssim 2\omega^2|| v||^2_{0,\Omega}+\big(\sum_{\gamma_{kj}} ||\lambda_h||^2_{0,\gamma_{kj}}\big)^\frac{1}{2} M^\frac{1}{2} \\
& \lesssim \omega^2h^{-1} M+\big(\sum_{\gamma_{kj}} ||\lambda_h||^2_{0,\gamma_{kj}}\big)^\frac{1}{2} M^\frac{1}{2}.
\end{split}
\end{equation}
Substituting (\ref{5.2new00}) of Lemma \ref{lamda} into (\ref{b4}), yields
\begin{equation*}
\||v\||^2\lesssim \omega^2 h^{-1} M + (h^{-1}pq\||v\||^2)^{{1\over 2}} M^\frac{1}{2}.
\end{equation*}
From the above inequality, we can deduce that
\begin{equation*}
\||v\||\lesssim \omega^2 h^{-1} M,
\end{equation*}
which gives the desired result (\ref{4.inequality.new10})
$\hfill \Box$

%%%%%%%%%%%%%%%%%%%%%%%%%%%%%%%%%%%%%%%%%%%%%%%%%%%%%%%%%%%%%%%%%%%%%%%%%%%%%%%%%%%%%%%
Now we can easily prove Theorem \ref{well} by Lemma \ref{lamda} and {\bf Proposition 5.1}.
\paragraph*{{\bf Proof of Theorem \ref{well}}.}
For $\lambda_h\in W^q_h(\gamma)$, let $u^{(1)}_{h,k}(\lambda_h)$ be the function defined in Subsection \ref{uh}.
From the definition of $u^{(1)}_{h,k}(\lambda_h)$, we have
\begin{equation*}
a^{(k)}(u_{h,k}^{(1)}(\lambda_h),\overline{w}_{h})
=\langle \pm\lambda_h, \overline{w}_{h}\rangle_{\partial\Omega_k \backslash\partial\Omega},~~k=1,\cdots,N;~~\forall\,w_h\in V_h^p(\mathcal{T}_h).
\end{equation*}
Namely, $u_{h}^{(1)}(\lambda_h)$ satisfies (\ref{bo}). It follows by Lemma \ref{lamda} that
\begin{equation} \label{new}
\sum_{\gamma_{kj}} ||\lambda_h||^2_{0,\gamma_{kj}}\lesssim h^{-1}pq\||u^{(1)}_h(\lambda_h)\||^2.
\end{equation}
Obviously,  $u_{h}^{(1)}(\lambda_h)$ satisfies (\ref{bb}) too. It follows by {\bf Proposition 5.1} that
\begin{equation*}
\||u^{(1)}_h(\lambda_h)\||^2
\lesssim \omega^2h^{-1} \sum_{\gamma_{kj}} ||[u^{(1)}_h(\lambda_h)]||^2_{0,\gamma_{kj}}.
\end{equation*}
This, together with (\ref{new}), leads to
\begin{equation*}
\sum_{\gamma_{kj}} ||\lambda_h||^2_{0,\gamma_{kj}}
\lesssim \omega^2h^{-2}pq\sum_{\gamma_{kj}} ||[u^{(1)}_h(\lambda_h)]||^2_{0,\gamma_{kj}}.
\end{equation*}
Thus
\begin{equation*}
s_h(\lambda_h,\lambda_h)
= \sum_{\gamma_{kj}} ||[u^{(1)}_h(\lambda_h)]||^2_{0,\gamma_{kj}}
\gtrsim \omega^{-2}h^{2}p^{-1}q^{-1} \sum_{\gamma_{kj}} ||\lambda_h||^2_{0,\gamma_{kj}}.
\end{equation*}
%where $\delta=1$ when $p\geq q+2$ ($q\geq 1$) or $\delta=0$ if $p=2q$ ($q\geq 2$).
$\hfill \Box$
\begin{remark} Remark 5.1 tells us that, when $q\geq 2$ and $p=2q$, a slightly better result than (\ref{4.new1}) can be built
$$ s_h(\lambda_h,\lambda_h)
= \sum_{\gamma_{kj}} ||[u^{(1)}_h(\lambda_h)]||^2_{0,\gamma_{kj}}
\gtrsim \omega^{-2}h^{2}p^{-1} \sum_{\gamma_{kj}} ||\lambda_h||^2_{0,\gamma_{kj}}.$$
\end{remark}

\subsection{Analysis on the error estimates}

In order to prove Theorem 4.3, we need more auxiliary results. We will decompose the error $u-u_h$ into three parts, where the first part and the second part have some particular property and the third
part is a $p$-order finite element function. The third part can be estimated by {\bf Proposition 5.1}, but the estimates of the first part and the second part are more technical, which depend
on a key auxiliary result (Lemma \ref{aaa}).

We first build the key auxiliary result mentioned above. For an element $\Omega_k$ and a function $v\in H^1(\Omega_k)$, we use the notation in this subsection
\begin{equation}  \label{a2}
F_k(v)=||\nabla v||^2_{0,\Omega_k}-(\kappa^2v,v)_{\Omega_k} \pm i\rho ||v||^2_{0,\partial\Omega_k \backslash\partial\Omega}+i\langle\kappa v,v\rangle_{\partial\Omega_k \cap\partial\Omega}.
\end{equation}
It is clear that $F_k(v)=a^{(k)}(v,v)$.

\begin{lemma} \label{aaa} Let Assumption 2 and Assumption 3 be satisfied. For one element $\Omega_k$, assume that $v\in H^1(\Omega_k)$ has the property
\begin{equation} \label{a1}
-(\kappa^2v, 1)_{\Omega_k} \pm i\rho\langle v, 1\rangle_{\partial\Omega_k \backslash\partial\Omega}+i\langle\kappa v, 1\rangle_{\partial\Omega_k \cap\partial\Omega} =0.
\end{equation}
Then
\begin{equation}
\|\nabla v\|^2_{0,\Omega_k}+\omega\|v\|^2_{0,\Omega_k}\leq C |F_k(v)|.\label{5.2new03}
\end{equation}
\end{lemma}

\paragraph*{Proof.} Taking the module to (\ref{a2}) leads to
\begin{equation} \label{5.2new01}
\big| ||\nabla v||^2_{0,\Omega_k}-(\kappa^2 v,v)_{\Omega_k}\big|\leq |F_k(v)|
\end{equation}
and
\begin{equation} \label{5.2new02}
\big|\pm\rho ||v||^2_{0,\partial\Omega_k \backslash\partial\Omega}+\langle\kappa v,v\rangle_{\partial\Omega_k \cap\partial\Omega}\big|
\leq |F_k(v)|.
\end{equation}

We first assume that $\partial\Omega_k\cap\partial\Omega \ne \emptyset$. It follows, by (\ref{5.2new02}) and the trace inequality, that
\begin{equation} \label{a5}
\begin{split}
\omega||v||^2_{0,\partial\Omega_k\cap\partial\Omega}
& %\lesssim ||\kappa^\frac{1}{2}v||^2_{0,\partial\Omega_k\cap\partial\Omega}
\lesssim |F_k(v)|+\rho ||v||^2_{0,\partial\Omega_k \backslash\partial\Omega} \\
&\lesssim |F_k(v)|+\rho h ||\nabla v||^2_{0,\Omega_k}+\rho h^{-1} ||v||^2_{0,\Omega_k}.
\end{split}
\end{equation}
Using Poincar$\acute{\text{e}}$ inequality and (\ref{a5}), yields
\begin{equation*}
\begin{split}
\omega^2||v||^2_{0,\Omega_k}
& \lesssim \omega^2h^2||\nabla v||^2_{0,\Omega_k} +\omega^2h||v||^2_{0,\partial\Omega_k\cap\partial\Omega} \\
& \lesssim \omega^2h^2||\nabla v||^2_{0,\Omega_k} +\omega h|F_k(v)|+\rho \omega h^2 ||\nabla v||^2_{0,\Omega_k}+\rho \omega||v||^2_{0,\Omega_k},
\end{split}
\end{equation*}
which implies that
\begin{equation*}
(1-\rho C\omega^{-1})\omega^2||v||^2_{0,\Omega_k}
\lesssim (\omega^2h^2+\rho\omega h^2)||\nabla v||^2_{0,\Omega_k} +\omega h|F_k(v)|.
\end{equation*}
Then, from Assumption 3, we have
\begin{equation} \label{fri}
\omega^2||v||^2_{0,\Omega_k}
\lesssim \omega^2h^2 ||\nabla v||^2_{0,\Omega_k} +\omega h|F_k(v)|.
\end{equation}
On the other hand, by (\ref{5.2new01}) and (\ref{fri}), we deduce that
\begin{equation*}
||\nabla v||^2_{0,\Omega_k}
\lesssim\omega^2||v||^2_{0,\Omega_k} + |F_k(v)|
\lesssim \omega^2h^2 ||\nabla v||^2_{0,\Omega_k}+(1+\omega h)|F_k(v)|,
\end{equation*}
which gives
\begin{equation*}
(1-C\omega^2h^2) ||\nabla v||^2_{0,\Omega_k}
\lesssim |F_k(v)|.
\end{equation*}
This, together with (\ref{fri}), leads to
\begin{equation*}
(1-C\omega^2h^2) \omega^2||v||^2_{0,\Omega_k}
\lesssim (\omega^2h^2+\omega h(1-C\omega^2h^2))|F_k(v)|.
\end{equation*}
Using Assumption 2, the above two inequalities give (\ref{5.2new03}) when $\partial\Omega_k\cap\partial\Omega \ne \emptyset$.

In the following we assume that $\partial\Omega_k\cap\partial\Omega = \emptyset$. It follows by (\ref{a1}) that
\begin{equation*}
\omega^2|\Omega_k| |\gamma_{\Omega_k}(v)|
\lesssim \rho |\langle v, 1\rangle_{\partial\Omega_k\backslash\partial\Omega}|
\lesssim \rho |\partial\Omega_k|^\frac{1}{2} ||v||_{0, \partial\Omega_k}
\end{equation*}
where $\gamma_{\Omega_k}(v) = \frac{1}{|\Omega_k|}\int_{\Omega_k} v \, dx$. Thus
\begin{equation*}
\omega^4|\gamma_{\Omega_k}(v)|^2
\lesssim \rho^2 |\partial\Omega_k|\cdot |\Omega_k|^{-2} ||v||^2_{0, \partial\Omega_k}.
\end{equation*}
This, together with the trace inequality (or $\varepsilon-$inequality), leads to
\begin{equation} \label{eps}
\begin{split}
\omega^4||\gamma_{\Omega_k}(v)||^2_{0,\Omega_k}
&\lesssim \rho^2 |\partial\Omega_k|\cdot|\Omega_k|^{-1} ||v||^2_{0, \partial\Omega_k} \\
& \lesssim \rho^2 h^{-1} (h||\nabla v||^2_{0,\Omega_k} +h^{-1} ||v||^2_{0,\Omega_k}) \\
& = \rho^2||\nabla v||^2_{0,\Omega_k} +\rho^2h^{-2} ||v||^2_{0,\Omega_k}.
\end{split}
\end{equation}
% Here we have used the assumption $\omega h\leq C_0$ when we use $\varepsilon-$inequality.
Using Friedrichs' inequality and (\ref{eps}), we deduce that
\begin{equation*}
\begin{split}
\omega^4||v||^2_{0,\Omega_k}
& \leq \omega^4||v-\gamma_{\Omega_k}(v)||^2_{0,\Omega_k} +\omega^4||\gamma_{\Omega_k}(v)||^2_{0,\Omega_k} \\
& \lesssim \omega^4h^2 ||\nabla v||^2_{0,\Omega_k}+\rho^2||\nabla v||^2_{0,\Omega_k} +\rho^2h^{-2} ||v||^2_{0,\Omega_k}.
\end{split}
\end{equation*}
So we get
\begin{equation*}
(1-\rho^2\omega^{-4}h^{-2})\omega^2 ||v||^2_{0,\Omega_k}
\lesssim (\omega^2h^2+\rho^2\omega^{-2}) ||\nabla v||^2_{0,\Omega_k}.
\end{equation*}
Thus, by Assumption 3
%the condition $\rho \ll \omega^2h$
, we have
\begin{equation*}
\omega^2||v||^2_{0,\Omega_k}
\lesssim \omega^2h^2 ||\nabla v||^2_{0,\Omega_k}.
\end{equation*}
In addition, combining (\ref{5.2new01}) with the above inequality, we get
\begin{equation*}
||\nabla v||^2_{0,\Omega_k}
\lesssim \omega^2||v||^2_{0,\Omega_k} + |F_k(v)|
\lesssim \omega^2h^2 ||\nabla v||^2_{0,\Omega_k}+|F_k(v)|.
\end{equation*}
Therefore we obtain (if $C\omega^2h^2<1$)
\begin{equation*}
(1-C\omega^2h^2) ||\nabla v||^2_{0,\Omega_k}
\lesssim |F_k(v)|
\end{equation*}
and
\begin{equation*}
(1-C\omega^2h^2) \omega^2||v||^2_{0,\Omega_k}
\lesssim \omega^2h^2|F_k(v)|.
\end{equation*}
Using Assumption 2 again, the above two inequalities give (\ref{5.2new03}) for the case that $\partial\Omega_k\cap\partial\Omega = \emptyset$.
$\hfill \Box$

Next we use Lemma \ref{aaa} to build two new auxiliary results, which involve an auxiliary function $\hat{u}_h(\lambda)$ for $\lambda\in W(\gamma)$.
For each element $\Omega_k$, let $\hat{u}_{h,k}(\lambda)\in V_h^p(\Omega_k)$ be determined by the variational problem
\begin{equation} \label{inter}
a^{(k)}(\hat{u}_{h,k}(\lambda),\overline{v}_{h})
=L^{(k)}(\overline{v}_{h})+\langle \pm\lambda, \overline{v}_{h}\rangle_{\partial\Omega_k \backslash\partial\Omega}, \quad\forall\,v_{h}\in V_h^p(\Omega_k).
\end{equation}
Then define
$\hat{u}_h(\lambda)\in V_h^p(\mathcal{T}_h)$ such that $\hat{u}_h(\lambda)|_{\Omega_k}=\hat{u}_{h,k}(\lambda)$ ($k=1,\cdots,N$).

\begin{lemma} \label{step1} Assume that $u\in H^{r+1}({\mathcal T}_h)$ with $1\leq r\leq p$.
Let Assumption 2 and Assumption 3 be satisfied. Then
\begin{equation*}
|u(\lambda)-\hat{u}_h(\lambda)|_{1,\Omega}
\leq C h^rp^{-r}|u|_{r+1,\Omega}
\end{equation*}
and
\begin{equation*}
||u(\lambda)-\hat{u}_h(\lambda)||_{0,\Omega}
\leq C \omega^{-1} h^rp^{-r}|u|_{r+1,\Omega}.
\end{equation*}

\end{lemma}

\paragraph*{Proof.}
Let $\varepsilon_u = u(\lambda)-\hat{u}_h(\lambda)$ and $\gamma_{\Omega_k}(\varepsilon_u) = \frac{1}{|\Omega_k|}\int_{\Omega_k}\varepsilon_u dx$.
Choosing $v=v_h$ in (\ref{viar1}) and taking the difference between (\ref{viar1}) and (\ref{inter}), we can get
\begin{equation} \label{1}
a^{(k)}(\varepsilon_u, \overline{v}_h)
= 0 \quad\quad  \forall\,v_h\in V_h^p(\Omega_k), ~k=1,2,\dots,N.
\end{equation}
Set $v_h=1$ in the above inequality, we have
\begin{equation*}
-(\kappa^2\varepsilon_u, 1)_{\Omega_k} \pm i\rho\langle \varepsilon_u, 1\rangle_{\partial\Omega_k \backslash\partial\Omega}+i\langle\kappa\varepsilon_u, 1\rangle_{\partial\Omega_k \cap\partial\Omega} =0.
%-\kappa^2|\Omega_k|\gamma_{\Omega_k}(\varepsilon_u) \pm i\rho\int_{\partial\Omega_k \backslash\partial\Omega}\varepsilon_u ds +i\kappa\int_{\partial\Omega_k \cap\partial\Omega}\varepsilon_u ds = 0
\end{equation*}
Then $\varepsilon_u$ has the property (\ref{a1}), so by Lemma \ref{aaa} we obtain
\begin{equation}
||\nabla \varepsilon_u||^2_{0,\Omega_k}+
\omega^2||\varepsilon_u||^2_{0,\Omega_k}
\lesssim |F_k(\varepsilon_u)|.\label{5.3new00}
\end{equation}
It suffices to estimate $|F_k(\varepsilon_u)|$.

Let $Q_h^p: L^2(\Omega_k) \mapsto V_h^p(\Omega_k)$ denote the standard $L^2$ projection operator. It is clear that
\begin{equation*}
\varepsilon_u = (I-Q_h^p)u(\lambda)+Q_h^p u(\lambda)-\hat{u}_h(\lambda).
\end{equation*}
Choosing $v_h=Q_h^p u(\lambda)-\hat{u}_h(\lambda)$ in (\ref{1}), we have
\begin{equation*}
a^{(k)}(\varepsilon_u, \overline{Q_h^p u(\lambda)-\hat{u}_h(\lambda)})
= 0,\quad \quad ~k=1,2,\dots,N.
\end{equation*}
Thus
$$ F_k(\varepsilon_u)=a^{(k)}(\varepsilon_u, \overline{\varepsilon_u})=a^{(k)}(\varepsilon_u, \overline{(I-Q_h^p)u(\lambda)}). $$
Using the definition of $a^{(k)}(\cdot,\cdot)$, Cauchy inequality and the assumption $\rho\leq \omega$, we further get
\begin{eqnarray}
|F_k(\varepsilon_u)|&\lesssim&|\varepsilon_u|_{1,\Omega_k}\cdot|(I-Q_h^p)u(\lambda)|_{1,\Omega_k}+\omega^2\|\varepsilon_u\|_{0,\Omega_k}\cdot\|(I-Q_h^p)u(\lambda)\|_{0,\Omega_k}\cr
&+&\omega\|\varepsilon_u\|_{0,\partial\Omega_k}\cdot\|(I-Q_h^p)u(\lambda)\|_{0,\partial\Omega_k}
\label{5.3new02}
\end{eqnarray}
By the trace inequality, we have
$$ \|\varepsilon_u\|_{0,\partial\Omega_k}\lesssim h^{{1\over 2}}|\varepsilon_u|_{1,\Omega_k}+h^{-{1\over 2}}\|\varepsilon_u\|_{0,\Omega_k} $$
and
$$ \|(I-Q_h^p)u(\lambda)\|_{0,\partial\Omega_k}\lesssim h^{{1\over 2}}|(I-Q_h^p)u(\lambda)|_{1,\Omega_k}+h^{-{1\over 2}}\|(I-Q_h^p)u(\lambda)\|_{0,\Omega_k}. $$
Substituting the above two inequalities into (\ref{5.3new02}) and using Schwarz inequality, yields
\begin{eqnarray}
&&|F_k(\varepsilon_u)|\lesssim(|\varepsilon_u|^2_{1,\Omega_k}+\omega^2\|\varepsilon_u\|^2_{0,\Omega_k})^{{1\over 2}}\cr
&&\cdot\big((1+\omega^2h^2)|(I-Q_h^p)u(\lambda)|^2_{1,\Omega_k}+(\omega^2+h^{-2})\|(I-Q_h^p)u(\lambda)\|^2_{0,\Omega_k}\big)^{{1\over 2}}.\label{5.3new03}
\end{eqnarray}

Using the approximation of the projection operator (refer to \cite{L2}), we have
\begin{equation*}
||\nabla(I-Q_h^p)u(\lambda)||_{0,\Omega_k}
\lesssim h^{r} p^{-r} |u(\lambda)|_{r+1,\Omega_k}
\end{equation*}
and
\begin{equation*}
||(I-Q_h^p)u(\lambda)||_{0,\Omega_k}
\lesssim h^{r+1} p^{-(r+1)} |u(\lambda)|_{r+1,\Omega_k}.
\end{equation*}
Plugging these in (\ref{5.3new03}), leads to
\begin{equation*}
|F_k(\varepsilon_u)|\lesssim(|\varepsilon_u|^2_{1,\Omega_k}+\omega^2\|\varepsilon_u\|^2_{0,\Omega_k})^{{1\over 2}}\cdot
\big((1+\omega^2 h^2+{1\over p^2}+{\omega^2 h^2\over p^2})h^{2r} p^{-2r} |u(\lambda)|^2_{r+1,\Omega_k}\big)^{{1\over 2}}.
\end{equation*}
Furthermore, by Assumption 2 we get
$$ |F_k(\varepsilon_u)| \lesssim (|\varepsilon_u|^2_{1,\Omega_k}+\omega^2\|\varepsilon_u\|^2_{0,\Omega_k})^{{1\over 2}}\cdot h^{r} p^{-r} |u(\lambda)|_{r+1,\Omega_k}, $$
which, together with (\ref{5.3new00}), gives
\begin{equation} \label{55}
|\varepsilon_u|^2_{1,\Omega_k}+\omega^2\|\varepsilon_u\|^2_{0,\Omega_k}\lesssim h^{2r} p^{-2r} |u(\lambda)|^2_{r+1,\Omega_k}.
\end{equation}
Then the estimates in this lemma can be obtained by summing (\ref{55}) over $k$.
$\hfill \Box$

Let $Q_h^q: L^2(\gamma)\rightarrow W_h^q(\gamma)$ denote the $L^2$ projector, and let $u_h(Q_h^q\lambda)\in V_h^p(\Omega_k)$ be defined as in (\ref{viar2}), by choosing $\lambda_h=Q_h^q\lambda$.
\begin{lemma} \label{step2} Suppose that $q\geq 1$, $p\geq q+2$ and $\lambda\in W^{r-{1\over 2}}(\gamma)$ with $1\leq r\leq q+1$.
Let Assumption 1-Assumption 4 be satisfied.
%Assume that $\omega h\leq C_0$ for a positive constant $C_0$
Then
\begin{equation*}
|\hat{u}_h(\lambda)-u_h(Q_h^q\lambda)|_{1,\Omega}
\leq C h^rq^{-r} |\lambda|_{r-\frac{1}{2},\gamma}
\end{equation*}
and
\begin{equation*}
||\hat{u}_h(\lambda)-u_h(Q_h^q\lambda)||_{0,\Omega}
\leq C \omega^{-1}h^rq^{-r} |\lambda|_{r-\frac{1}{2},\gamma}.
\end{equation*}

\end{lemma}

\paragraph*{Proof.} Set $\tilde{\varepsilon}_u = \hat{u}_h(\lambda)-u_h(Q_h^q\lambda)$ and $\varepsilon_{\lambda} = \lambda-Q_h^q\lambda$. From (\ref{viar2}) and (\ref{inter}), we deduce that
\begin{equation} \label{7}
\begin{split}
a^{(k)}(\tilde{\varepsilon}_u, \overline{v}_h)
= \langle \pm\varepsilon_{\lambda}, \overline{v}_h\rangle_{\partial\Omega_k \backslash\partial\Omega}, \quad \forall\,v_h\in V_h^p(\Omega_k), ~ k=1,2,\dots,N.
\end{split}
\end{equation}
In particular, choosing $v_h=1$ in (\ref{7}), we have
\begin{equation*}
-(\kappa^2\tilde{\varepsilon}_u, 1)_{\Omega_k} \pm i\rho\langle \tilde{\varepsilon}_u, 1\rangle_{\partial\Omega_k \backslash\partial\Omega}+i\langle\kappa\tilde{\varepsilon}_u, 1\rangle_{\partial\Omega_k \cap\partial\Omega}
= \int_{\partial\Omega_k \backslash\partial\Omega} \pm\varepsilon_{\lambda} ds
= 0.
\end{equation*}
Then, by applying Lemma \ref{aaa} to $\tilde{\varepsilon}_u$, we get
\begin{equation}
||\nabla \tilde{\varepsilon}_u||^2_{0,\Omega_k}+
\omega^2||\tilde{\varepsilon}_u||^2_{0,\Omega_k}
\lesssim |F_k(\tilde{\varepsilon}_u)|.\label{5.3new04}
\end{equation}

On the other hand, set $v_h = \tilde{\varepsilon}_u$ in (\ref{7}), then we get
\begin{equation*}
\begin{split}
&F_k(\tilde{\varepsilon}_u)\doteq||\nabla \tilde{\varepsilon}_u||^2_{0,\Omega_k}-(\kappa^2\tilde{\varepsilon}_u,\tilde{\varepsilon}_u)_{\Omega_k} \pm i\rho ||\tilde{\varepsilon}_u||^2_{0,\partial\Omega_k \backslash\partial\Omega}+i\langle\kappa\tilde{\varepsilon}_u,\tilde{\varepsilon}_u\rangle_{\partial\Omega_k \cap\partial\Omega} \\
& = \langle \pm\varepsilon_{\lambda}, \tilde{\varepsilon}_u \rangle_{\partial\Omega_k\backslash\partial\Omega}
= \langle \pm\varepsilon_{\lambda}, \tilde{\varepsilon}_u-\gamma_{\Omega_k}(\tilde{\varepsilon}_u) \rangle_{\partial\Omega_k\backslash\partial\Omega}.
\end{split}
\end{equation*}
Here we have used the fact that $\gamma_{\Omega_k}(\tilde{\varepsilon}_u) = \frac{1}{|\Omega_k|}\int_{\Omega_k}\tilde{\varepsilon}_u dx$ is a constant. It is easy to see that
\begin{equation*}
|F_k(\tilde{\varepsilon}_u)|\leq ||\pm\varepsilon_{\lambda}||_{-\frac{1}{2},\partial\Omega_k\backslash\partial\Omega} ||\tilde{\varepsilon}_u -\gamma_{\Omega_k}(\tilde{\varepsilon}_u)||_{\frac{1}{2},\partial\Omega_k}
\lesssim ||\pm\varepsilon_{\lambda}||_{-\frac{1}{2},\partial\Omega_k\backslash\partial\Omega} ||\nabla\tilde{\varepsilon}_u||_{0,\Omega_k}.
\end{equation*}
Substituting this into (\ref{5.3new04}), yields
\begin{equation}
||\nabla \tilde{\varepsilon}_u||^2_{0,\Omega_k}+\omega^2||\tilde{\varepsilon}_u||^2_{0,\Omega_k}
\lesssim ||\pm\varepsilon_{\lambda}||^2_{-\frac{1}{2},\partial\Omega_k\backslash\partial\Omega}.\label{5.inequality-new11}
\end{equation}

By using the approximation of the projection operator $Q_h^q$, we have
\begin{equation*}
||\pm\varepsilon_{\lambda}||_{-\frac{1}{2},\partial\Omega_k\backslash\partial\Omega}
\lesssim h^rq^{-r} |\lambda|_{r-\frac{1}{2},\partial\Omega_k\backslash\partial\Omega}.
\end{equation*}
Combining this with (\ref{5.inequality-new11}), leads to
\begin{equation*} \label{4}
|\nabla \tilde{\varepsilon}_u|^2_{1,\Omega_k}
\lesssim h^rq^{-r} |\lambda|_{r-\frac{1}{2},\partial\Omega_k\backslash\partial\Omega}
+\omega^2||\tilde{\varepsilon}_u||^2_{0,\Omega_k}
\lesssim h^{2r}q^{-2r} |\lambda|^2_{r-\frac{1}{2},\partial\Omega_k\backslash\partial\Omega}.
\end{equation*}
By summing the above inequalities over $k$, we obtain the desired estimates.
$\hfill \Box$

In the following we use {\bf Proposition 5.1} to build an estimate of $u_h(Q_h^q\lambda)-u_h(\lambda_h)$.

\begin{lemma} \label{step3}
Suppose that $q\geq 1$, $p\geq q+2$ and $u\in H^{r+1}({\mathcal T}_h)$ with $1\leq r\leq p$. Let Assumption 1-Assumption 4 be satisfied.
%$p\geq 1+c_0\log\omega$ for a constant $c_0$. Then, there exists a constant $C_0$ such that, when $\omega h\leq C_0$
Then
\begin{equation*}
|u_h(Q_h^q\lambda)-u_h(\lambda_h)|_{1,\Omega}
\leq Ch^{r-1} ( p^{-r}|u|_{r+1,\Omega} + q^{-r}|\lambda|_{r-\frac{1}{2},\gamma})
\end{equation*}
and
\begin{equation*}
%\sum_{k=1}^N \kappa_0||u_h(Q_h^q\lambda)-u_h(\lambda_h)||_{0,\Omega_k}
||u_h(Q_h^q\lambda)\!-\!u_h(\lambda_h)||_{0,\Omega}
\leq C\omega^{-1} h^{r-1} ( p^{-r}|u|_{r+1,\Omega} + q^{-r}|\lambda|_{r-\frac{1}{2},\gamma} ).
\end{equation*}

\end{lemma}

\paragraph*{Proof.}
Set $\varepsilon_h = u_h(Q_h^q\lambda)-u_h(\lambda_h)$ and let $\tilde{\varepsilon}_{\lambda} = Q_h^q\lambda-\lambda_h$. Notice that the function $u_h(Q_h^q\lambda)$
is defined by (\ref{viar2}) with $\lambda_h=Q_h^q\lambda$. Then, it follows by (\ref{viar2}) that
\begin{equation} \label{14}
a^{(k)}(\varepsilon_h, w)
= \langle \pm\tilde{\varepsilon}_{\lambda}, \bar{w}\rangle_{\partial\Omega_k\backslash\partial\Omega}, \quad \forall\,w\in V_h^p(\mathcal{T}_h), ~k=1,2,\dots,N.
\end{equation}
Setting $v = \varepsilon_h$ in {\bf Proposition 5.1}, we get
\begin{equation}
\sum_{k=1}^N||\nabla \varepsilon_h||^2_{0,\Omega}+\omega^2||\varepsilon_h||^2_{0,\Omega}=\||\varepsilon_h\||^2
\lesssim \omega^2h^{-1} \sum_{\gamma_{kj}} ||[\varepsilon_h]||^2_{0,\gamma_{kj}}.\label{5.inequality-new13}
\end{equation}

By the definition of $\lambda_h$, which corresponds to the minimal energy, we deduce that
\begin{equation}
\begin{split}
& \sum_{\gamma_{kj}}||[\varepsilon_h]||^2_{0,\gamma_{kj}}
\leq \sum_{\gamma_{kj}}||[u_h(Q_h^q\lambda)]||^2_{0,\gamma_{kj}} +\sum_{\gamma_{kj}}||[u_h(\lambda_h)]||^2_{0,\gamma_{kj}} \\
& \leq 2\sum_{\gamma_{kj}}||[u_h(Q_h^q\lambda)]||^2_{0,\gamma_{kj}}
=2\sum_{\gamma_{kj}}||[u_h(Q_h^q\lambda)-u(\lambda)]||^2_{0,\gamma_{kj}}.\label{5.inequality-new12}
\end{split}
\end{equation}
Here we have used the fact that the function $u(\lambda)$ has the zero jump across each edge $\gamma_{kj}$. Using the trace inequality ($\varepsilon$-inequality), yields
\begin{equation*}
\begin{split}
||[u_h(Q_h^q\lambda)-u(\lambda)]||^2_{0,\gamma_{kj}}
& \lesssim h||\nabla (u_h(Q_h^q\lambda)-u(\lambda))||^2_{0,\Omega_k} +h^{-1}||u_h(Q_h^q\lambda)-u(\lambda)||^2_{0,\Omega_k}\\
& \lesssim h(||\nabla (u_h(Q_h^q\lambda)-\hat{u}_h(\lambda))||^2_{0,\Omega_k}+||\nabla (\hat{u}_h(\lambda)-u(\lambda))||^2_{0,\Omega_k}) \\
& +h^{-1}(||u_h(Q_h^q\lambda)-\hat{u}_h(\lambda)||^2_{0,\Omega_k} +||\hat{u}_h(\lambda)-u(\lambda)||^2_{0,\Omega_k}).
\end{split}
\end{equation*}
Substituting this inequality into (\ref{5.inequality-new12}) and using Lemma \ref{step1} and Lemma \ref{step2}, yields
$$
\sum_{\gamma_{kj}}||[\varepsilon_h]||^2_{0,\gamma_{kj}}
%\lesssim \omega^{-1} ( h^{2r}p^{-2r}|u|^2_{r+1,\Omega} + h^{2r}q^{-2r}\|\lambda\|^2_{r-\frac{1}{2},\gamma} ).
\lesssim (h+h^{-1}\omega^{-2})( h^{2r}p^{-2r}|u|^2_{r+1,\Omega} + h^{2r}q^{-2r}|\lambda|^2_{r-\frac{1}{2},\gamma} ).
$$
This, together with (\ref{5.inequality-new13}), leads to
$$
\sum_{k=1}^N||\nabla \varepsilon_h||^2_{0,\Omega_k}+\omega^2||\varepsilon_h||^2_{0,\Omega}
\lesssim (\omega^2+h^{-2}) h^{2r}( p^{-2r}|u|^2_{r+1,\Omega} + q^{-2r}|\lambda|^2_{r-\frac{1}{2},\gamma} ).
$$
Then we immediately obtain the estimates in this lemma.
$\hfill \Box$

Now we can easily prove Theorem \ref{fina} by Lemma  5.7-Lemma 5.9.
\vspace {1mm}
\paragraph*{{\bf Proof of Theorem \ref{fina}}.} By the triangle inequality,  we have
\begin{equation*}
\begin{split}
& ||u(\lambda)-u_h(\lambda_h)||_{0,\Omega} \\
& \leq ||u(\lambda)-\hat{u}_h(\lambda)||_{0,\Omega} +||\hat{u}_h(\lambda)-u_h(Q_h^q\lambda)||_{0,\Omega} +||u_h(Q_h^q\lambda)-u_h(\lambda_h)||_{0,\Omega}
\end{split}
\end{equation*}
and
\begin{equation*}
\begin{split}
& |u(\lambda)-u_h(\lambda_h)|_{1,\Omega} \\
& \leq |u(\lambda)-\hat{u}_h(\lambda)|_{1,\Omega} +|\hat{u}_h(\lambda)-u_h(Q_h^q\lambda)|_{1,\Omega} +|u_h(Q_h^q\lambda)-u_h(\lambda_h)|_{1,\Omega}.
\end{split}
\end{equation*}
Then, by the estimates given in Lemma \ref{step1}, \ref{step2} and \ref{step3}, we obtain
\begin{equation*}
||u(\lambda)-u_h(\lambda_h)||_{0,\Omega}
\lesssim \omega^{-1} h^{r-1} ( p^{-r}|u|_{r+1,\Omega}+ q^{-r}|\lambda|_{r-\frac{1}{2},\gamma} )
\end{equation*}
and
\begin{equation*}
|u(\lambda)-u_h(\lambda_h)|_{1,\Omega}
\lesssim  h^{r-1} ( p^{-r}|u|_{r+1,\Omega}+ q^{-r}|\lambda|_{r-\frac{1}{2},\gamma} ).
\end{equation*}
$\hfill \Box$
\begin{remark} The proposed method can be extended directly to the case with other boundary conditions that can guarantee the well-posedness of the equations, provided
that the subproblems defined on the elements touching the boundary $\partial\Omega$ are imposed the corresponding boundary conditions (the analysis is simpler for the variational
problems with part Dirichlet boundary condition). The proposed discretization method can be also extended to three-dimensional problems, but the coarse subspace involved in the
construction of the preconditioner needs to be modified and the analysis is more difficult (for example,
the proofs of Theorem 4.1 and Lemma 5.3 needs to be changed).
\end{remark}

\section{Numerical experiments}
In this section we report some numerical results to illustrate that the proposed least squares method and domain decomposition preconditioner are efficient for Helmholtz
equations with large wave numbers.

In the discretization method described in Section 2, the parameter $\rho$ can be relatively arbitrarily positive number. We find that the different choices of $\rho$ does not
affect the accuracy of the resulting approximations provided that the value of $\rho$ is less than $1$. In this section, we simply choose $\rho=10^{-5}$ for numerical experiments.

For the considered example, the domain $\Omega$ is a rectangle so we adopt a uniform partition $\mathcal {T}_h$ for the domain $\Omega$ as follows:
$\Omega$ is divided into some small rectangles with a same size $h$, where $h$ denotes the length of the longest edge of the elements.

To measure the accuracy of the numerical solution $u_h$, we introduce the following relative $L^2$ error:
\begin{equation*}
\text{Err.}=\frac{||u_{ex}-u_h||_{L^2(\Omega)}}{||u_{ex}||_{L^2(\Omega)}}.
\end{equation*}

For a discretization method, when the value of $\omega h$ is fixed but $\omega$ increases ($h$ decreases), the relative $L^2$ error $\text{Err.}$ may obviously increase (if the number of basis
functions on each element does not increase).
This phenomenon is called ``wave number pollution". The efficiency of a discretization method for Helmholtz equations can be characterized by the degree of wave number pollution.
For convenience, we define a positive parameter $\delta$ to measure the {\it degree of wave number pollution} as follows: the parameter $\delta$ is the minimal positive number
such that, when $\omega$ increases and $h$ decreases to keep
the value of $\omega^{1+\delta} h$ being a constant, the relative $L^2$ error $\text{Err.}$ does not increase. If $\delta=0$, the discretization method has no ``wave number pollution".
For the standard linear finite element method, the existing results imply that $\delta=1$ (see \cite{fem}). A discretization method is ideal means that $\delta\ll 1$.
For concrete examples, it is difficult to exactly calculate such parameter $\delta$. Because of this, we want to give a similar definition of $\delta$, which can be explicitly calculated.

When $\omega$ increases from $\omega_1$ to $\omega_2$, the mesh size $h$ decreases from $h_1$ to $h_2$. We fix the value $\omega h$, i.e.,
$\omega_2h_2=\omega_1 h_1$. Let $\text{Err}_1$ and $\text{Err}_2$ denote the relative $L^2$ errors with $\omega=\omega_1$ ($h=h_1$) and $\omega=\omega_2$ ($h=h_2$), respectively.
Define $\delta>0$ by
$$ {\omega_2^{1+\delta}h_2\over \omega_1^{1+\delta}h_1}={ \text{Err}_2\over \text{Err}_1}. $$
It is easy to see that the parameter $\delta$ can be expressed as
$$ \delta = \frac{\ln(\text{Err}_2/\text{Err}_1)}{\ln(\omega_2/\omega_1)}+\bigg({\ln(h_1/h_2)\over \ln(\omega_2/\omega_1)}-1\bigg)=\frac{\ln(\text{Err}_2/\text{Err}_1)}{\ln(\omega_2/\omega_1)}
\quad\quad(\mbox{since}~~\omega_2h_2=\omega_1 h_1). $$

For a given $\omega$, we define the error order with respect to $h$ in the standard manner, namely,
$$ \text{order} = \frac{\ln(\text{Err}_2/\text{Err}_1)}{\ln(h_2/h_1)},$$
where $\text{Err}_1$ and $\text{Err}_2$ denote the relative $L^2$ errors corresponding to $h=h_1$ and $h=h_2$ ($p,q$ are fixed), respectively.

Throughout this section we can simply choose $p=q+2$ to avoid extra cost of calculation.

\subsection{Wave propagation in a duct with rigid walls}
In this subsection, we give some comparisons between the proposed method and the plane wave least squares (PWLS) method for a homogeneous Helmholtz equation with constant wave number.
For the comparisons, we recall the basic ideas of the PWLS method (see Subsection 2.3 in \cite{HuZ2016}). In the PWLS method, the solution space consists of plane wave basis functions that exactly satisfy the considered homogeneous Helmholtz equation,
and the variational formula is derived by a minimization problem with a quadratic subject functional defined by the jumps of function values and normal derivations across all the element interfaces.
Since the basis functions satisfy the considered Helmholtz equation, one needs not to introduce auxiliary unknowns on the element interfaces and so does not solve local Helmholtz equations on elements.

We consider the following model Helmholtz equation for the acoustic pressure $u$ (see \cite{Huttunen2009})
\begin{equation}
\left\{
\begin{array}{rr}
-\Delta u-\omega^2u=0 \quad \text{in} \quad \Omega,\\
{\partial u\over
\partial {\bf n}}+i\omega u=g \quad \text{on}
\quad \partial\Omega,
\end{array}
\right.
 \label{exam1}
\end{equation}
where $\Omega=[0,2]\times[0,1]$, and $g=({\partial \over
\partial {\bf n}}+i\omega)u_{ex}$. The analytic solution $u_{ex}$ of the problem can be obtained in the closed form as
$$u_{ex}(x,y)=\text{cos}(k\pi y)(A_1e^{-i\omega_x x}+A_2e^{i\omega_x
x}) $$ with $\omega_x=\sqrt{\omega^2-(k\pi)^2}$, and the coefficients
$A_1$ and $A_2$ satisfying the equation
\begin{equation}
\left( {\begin{array}{cc} \omega_x & -\omega_x \\
(\omega-\omega_x)e^{-2i\omega_x} & (\omega+\omega_x)e^{2i\omega_x}
\end{array} }
\right)
 \left ( {\begin{array}{c} A_1 \\ A_2
\end{array}}
\right ) = \left ( {\begin{array}{c}
-i   \\
0
\end{array}}
\right ).
\end{equation}
%Without loss of generality, we choose $k=12$ in the experiments for the example.
%The solution represents propagating modes and
%evanescent modes when the mode number $k$ is below the cut-off value
%$k< k_{\text{cut-off}}={\omega\over\pi}$.

Let ``NLS'' denote the novel least squares method proposed in this paper. Besides,
let $\hat{p}$ be the number of plane wave basis functions on every elements. For convenience, we use ``dof.'' to denote the number of degrees of freedom in the resulting algebraic systems
(which mean the system (\ref{axb}) for the NLS method).

In Table 1-Table 3, we compare the required numbers of degrees of freedom to achieve almost the same accuracies of the approximate solutions generated by the two methods.
\begin{center}
\tabcaption{}
\label{com1}
fixing $\omega = 20\pi$ and decreasing $h$ (setting $k=19$)
\vskip 0.1in
\begin{tabular}{|c|c|c|c|c|}
\hline
 &\multicolumn{2}{c|}{$\text{PWLS},~~\hat{p}=12$} &\multicolumn{2}{c|}{$\text{NLS},~~(q,p)\!=\!(3,\!5)$} \\
\hline
$h$ &dof. &Err. &dof. &Err. \\
\hline
$\frac{1}{28}$ &18816 &1.85e-4 &12208 &6.72e-5 \\
\hline
$\frac{1}{36}$ &31104 &2.75e-5 &20304 &1.61e-5 \\
\hline
$\frac{1}{44}$ &46464 &7.83e-6 &30448 &5.26e-6 \\
\hline
$\frac{1}{52}$ &64896 &2.81e-6 &42640 &2.11e-6 \\
\hline
\end{tabular}
\end{center}

\begin{center}
\tabcaption{}
\label{addcom}
fixing $\omega = 40\pi$ and decreasing $h$ (setting $k=19$)
\vskip 0.1in
\begin{tabular}{|c|c|c|c|c|}
\hline
 &\multicolumn{2}{c|}{$\text{PWLS},~~\hat{p}=12$} &\multicolumn{2}{c|}{$\text{NLS},~~(q,p)\!=\!(4,\!6)$} \\
\hline
$h$ &dof. &Err. &dof. &Err. \\
\hline
$\frac{1}{28}$ &18816 &1.24e-2 &15260 &4.57e-4 \\
\hline
$\frac{1}{36}$ &31104 &2.19e-3 &25380 &7.26e-5 \\
\hline
$\frac{1}{44}$ &46464 &2.67e-4 &38060 &1.80e-5 \\
\hline
$\frac{1}{52}$ &64896 &5.99e-5 &53300 &6.02e-6 \\
\hline
\end{tabular}
\end{center}

\begin{center}
\tabcaption{}
\label{com2}
fixing $\omega h = 5\pi/8$ and increasing $\omega$ (setting $k=12$)
\vskip 0.1in
\begin{tabular}{|c|c|c|c|c|c|}
\hline
& &\multicolumn{2}{c|}{$\text{PWLS},~~\hat{p}=12$} &\multicolumn{2}{c|}{$\text{NLS},~~(q,p)\!=\!(3,\!5)$} \\
\hline
$\omega$ &$h$ &dof. &Err. &dof. &Err. \\
\hline
$30\pi$ &$\frac{1}{48}$ &55296 &1.29e-4 &36288 &3.24e-5 \\
\hline
$35\pi$ &$\frac{1}{56}$ &75264 &3.43e-4 &49504 &3.56e-5 \\
\hline
$40\pi$ &$\frac{1}{64}$ &98304 &5.73e-4 &64768 &3.81e-5 \\
\hline
$45\pi$ &$\frac{1}{72}$ &124416 &7.85e-4 &82080 &4.00e-5 \\
\hline
\end{tabular}
\end{center}
\vskip 0.2in

It can be seen from the above datas that, for the new least squares method, less degrees of freedom in the solved algebraic system are enough to achieve almost the same accuracies (with the same choices of $\omega$ and $h$).
For the proposed method, a little extra cost is needed when solving all the local problems defined on the elements (in parallel). Besides, the system (\ref{axb}) has more complex structure than the one in
the PWLS method and so its preconditioner is more difficult to construct.
In summary, the proposed method is at least comparable to the plane wave method
even if the wave number is a constant (otherwise, the plane wave method may be unpractical).

\subsection{An example with variable wave numbers} In this subsection, we consider the following Helmholtz equations with variable wave numbers
\begin{equation}
\left\{
\begin{aligned}
& -\Delta u - \kappa^2 u = f
& \text{in}\ \Omega,\\
& \frac{\partial u}{\partial n} + i\kappa u = g
& \text{on}\ \partial\Omega,
\end{aligned}
\right.
\end{equation}
where $\Omega=[0,1]\times[0,1]$ and $\kappa=\frac{\omega}{c(x,y)}$. We define the velocity field $c({\bf x})$ as a smooth converging lens with a Gaussian profile at the center $(r_1,r_2) = (1/2,1/2)$ (refer to \cite{num})
\begin{equation}
c(x,y) = \frac{4}{3}\Big(1-\frac{1}{8}\exp\big(-32((x-r_1)^2+(y-r_2)^2)\big)\Big).
\end{equation}
The analytic solution of the problem is given by
\begin{equation}
u_{ex}(x,y) = c(x,y)\exp(i\omega xy).
\end{equation}

For this example, the standard plane wave methods are unpractical. In Table \ref{err1} and Table \ref{err2}, we list the accuracies of the approximate solutions generated by the proposed least squares method,
where the algebraic systems are solved in the exact manner.
\begin{center}
\tabcaption{}
\label{err1}
Degrees of wave number pollution: fixing $\omega h$ to be a constant and increasing $\omega$ (and decreasing $h$)
\vskip 0.1in
\begin{tabular}{|c|c|c|c|c|c|c|c|c|c|}
\hline
 &\multicolumn{3}{c|}{$\omega h = 1,~(q,p)=(2,4)$} &\multicolumn{3}{c|}{$\omega h = 2,~(q,p)=(3,5)$} &\multicolumn{3}{c|}{$\omega h = 2,~(q,p)=(4,6)$}  \\
\hline
$\omega$ &$h$ &$\text{Err.}$ &$\delta$ &$h$ &$\text{Err.}$ &$\delta$ &$h$ &$\text{Err.}$ &$\delta$  \\
\hline
64 &$\frac{1}{64}$ &3.079e-5 & &$\frac{1}{32}$ &2.8643e-5 & &$\frac{1}{32}$ &1.690e-6 & \\
\hline
128 &$\frac{1}{128}$ &3.132e-5 &0.024 &$\frac{1}{64}$ &2.913e-5 &0.035 &$\frac{1}{64}$ &1.731e-6 &0.034 \\
\hline
256 &$\frac{1}{256}$ &3.241e-5 &0.049 &$\frac{1}{128}$ &2.952e-5 &0.019 &$\frac{1}{128}$ &1.753e-6 &0.019 \\
\hline
512 &$\frac{1}{512}$ &3.318e-5 &0.034 &$\frac{1}{256}$ &2.981e-5 &0.014 &$\frac{1}{256}$ &1.770e-6 &0.014 \\
\hline
\end{tabular}
\end{center}

\begin{center}
\tabcaption{}
\label{err2}
Convergence orders of the approximations with respect to $h$: fixing $\omega=64$ and decreasing $h$
\vskip 0.1in
\begin{tabular}{|c|c|c|c|c|c|c|c|c|c|}
\hline
 &\multicolumn{3}{c|}{$(q,p)=(2,4)$} &\multicolumn{3}{c|}{$(q,p)=(3,5)$} &\multicolumn{3}{c|}{$(q,p)=(4,6)$}  \\
\hline
$\omega$ &$h$ &$\text{Err.}$ &order &$h$ &$\text{Err.}$ &order &$h$ &$\text{Err.}$ &order \\
\hline
64 &$\frac{1}{32}$ &4.484e-4 & &$\frac{1}{16}$ &1.186e-3 & &$\frac{1}{16}$ &1.381e-4 & \\
\hline
64 &$\frac{1}{64}$ &3.079e-5 &3.864 &$\frac{1}{32}$ &2.843e-5 &5.383 &$\frac{1}{32}$ &1.690e-6 &6.352 \\
\hline
64 &$\frac{1}{128}$ &2.016e-6 &3.933 &$\frac{1}{64}$ &7.390e-7 &5.266 &$\frac{1}{64}$ &2.639e-8 &6.001 \\
\hline
64 &$\frac{1}{256}$ &1.282e-7 &3.975 &$\frac{1}{128}$ &2.174e-8 &5.087 &$\frac{1}{128}$ &4.135e-10 &5.996 \\
\hline
64 &$\frac{1}{512}$ &8.068e-9 &3.990 &$\frac{1}{256}$ &6.710e-10 &5.018 &$\frac{1}{256}$ &6.731e-12 &5.941\\
\hline
\end{tabular}
\end{center}
\vskip 0.1in

From the above two tables, we can see that the approximate solutions generated by the proposed method indeed have high accuracies and have little ``wave number pollution".

Since the resulting stiffness matrix is Hermitian positive definite, we can solve the system by the CG method and the PCG method with the preconditioner constructed in Section \ref{pre}.
As usual we choose $d\approx\sqrt{h}$ as the subdomain size in this preconditioner to guarantee the loading balance.
The stopping criterion in the iterative algorithms is that the relative $L^2$-norm of the residual of the iterative approximation satisfies $\epsilon < 1.0e-6 \ $.

Moreover, let $N_{iter}^{CG}$ represent the iteration count for solving the algebraic system by CG method and $N_{iter}^{PCG}$ represent the iteration count for solving the algebraic system by PCG method with the DD preconditioner. When the wave number $\omega$ increases (and the mesh size $h$ decreases), the iteration count $N_{iter}$ (represent $N_{iter}^{CG}$ or $N_{iter}^{PCG}$) also increases. In order to describe the growth rate of the iteration count $N_{iter}$ with respect to the wave number $\omega$, we introduce a new notation $\rho^{iter}$. Let $\omega_1$ and $\omega_2$ be two wave numbers, and let $N_{iter}^{(1)}$ and $N_{iter}^{(2)}$ denote the corresponding iteration counts, respectively. Then we define the positive number $\rho^{iter}$ by
$$ \frac{N_{iter}^{(2)}}{N_{iter}^{(1)}} = \big(\frac{\omega_2}{\omega_1}\big)^{\rho^{iter}}.$$

For example, when $\rho^{iter}=1$, the growth is linear; if $\rho^{iter}\to 0^+$, then the preconditioner possesses the optimal convergence. For a preconditioner, the positive number $\rho^{iter}$ defined above is called ``relative growth rate'' of the iteration count. Of course, we hope that the relative growth rate $\rho^{iter}$ is small.

In Table \ref{tab1}, Table \ref{tab2} and Table \ref{addtab}, we compare the iteration counts and its ``relative growth rate'' for the CG method and PCG method with the DD preconditioner constructed in Section 3.
\begin{center}
\tabcaption{}
\label{tab1}
Effectiveness of the preconditioner: the case with $\omega h\approx 1$ and $(q,p)=(2,4)$
\vskip 0.1in
\begin{tabular}{|c|c|c|c|c|c|c|c|}
\hline
$\omega$ &$h$ &$d$ &$N_{iter}^{CG}$ &$\rho^{iter}_{CG}$ &$N_{iter}^{PCG}$ &$\rho^{iter}_{PCG}$ &$\text{Err.}$ \\%&$\text{Err.}$
\hline
20$\pi$ & $\frac{1}{64}$ & $\frac{1}{8}$ &1556 & &105 &  &2.8825e-5  \\%&2.9833e-5
\hline
40$\pi$ & $\frac{1}{121}$ & $\frac{1}{11}$ &2504 &0.6864 &139 &0.4047 &3.8198e-5  \\%&3.9091e-5
\hline
80$\pi$ & $\frac{1}{256}$ & $\frac{1}{16}$ &5158 &1.0426 &191 &0.4585 &3.1632e-5 \\%&4.8094e-5
\hline
160$\pi$ & $\frac{1}{484}$ & $\frac{1}{22}$ &9643 &0.9027 &251 &0.3941 &4.2254e-5  \\
\hline
\end{tabular}
\end{center}

\begin{center}
\tabcaption{}
\label{tab2}
Effectiveness of the preconditioner: the case with $\omega h\approx 2$ and $(q,p)=(3,5)$
\vskip 0.1in
\begin{tabular}{|c|c|c|c|c|c|c|c|}
\hline
$\omega$ &$h$ &$d$ &$N_{iter}^{CG}$ &$\rho^{iter}_{CG}$ &$N_{iter}^{PCG}$ &$\rho^{iter}_{PCG}$ &$\text{Err.}$ \\%&$\text{Err.}$
\hline
20$\pi$ & $\frac{1}{36}$ & $\frac{1}{6}$ &451 & &77 & &1.4064e-5  \\%&1.4928e-5
\hline
40$\pi$ & $\frac{1}{64}$ & $\frac{1}{8}$ &602 &0.4166 &104 &0.4337 &2.7665e-5  \\%&3.0537e-5
\hline
80$\pi$ & $\frac{1}{121}$ & $\frac{1}{11}$ &967 &0.6838 &144 &0.4695 &3.9366e-5 \\%&4.1386e-5
\hline
160$\pi$ & $\frac{1}{256}$ & $\frac{1}{16}$ &1714 &0.8258 &190 &0.3999 &3.6951e-5  \\
\hline
\end{tabular}
\end{center}
\vskip 0.1in

\begin{center}
\tabcaption{}
\label{addtab}
Effectiveness of the preconditioner: the case with $\omega h\approx 2$ and $(q,p)=(4,6)$
\vskip 0.1in
\begin{tabular}{|c|c|c|c|c|c|c|c|}
\hline
$\omega$ &$h$ &$d$ &$N_{iter}^{CG}$ &$\rho^{iter}_{CG}$ &$N_{iter}^{PCG}$ &$\rho^{iter}_{PCG}$ &$\text{Err.}$ \\%&$\text{Err.}$
\hline
20$\pi$ & $\frac{1}{36}$ & $\frac{1}{6}$ &489 & &78 & &7.3846e-7  \\%&1.4928e-5
\hline
40$\pi$ & $\frac{1}{64}$ & $\frac{1}{8}$ &646 &0.4017 &106 &0.4425 &1.5461e-6  \\%&3.0537e-5
\hline
80$\pi$ & $\frac{1}{121}$ & $\frac{1}{11}$ &1002 &0.6333 &144 &0.4420 &2.2087e-6 \\%&4.1386e-5
\hline
160$\pi$ & $\frac{1}{256}$ & $\frac{1}{16}$ &1758 &0.8111 &191 &0.4075 &1.5775e-6  \\
\hline
\end{tabular}
\end{center}
\vskip 0.1in

The above data indicate that the proposed preconditioner is very efficient and the iteration counts of the corresponding PCG method has small relative
growth rate when the wave number increases.

\bibliographystyle{siamplain}

\end{document}